\documentclass[12pt]{amsart}

\usepackage{amsmath,amssymb,amsthm,mathtools}
\usepackage[margin=1in]{geometry}
\usepackage{enumitem}
\usepackage{microtype}
\usepackage{tikz}
\usepackage[hidelinks]{hyperref}

\theoremstyle{plain}
\newtheorem{theorem}{Theorem}[section]
\newtheorem{proposition}[theorem]{Proposition}
\newtheorem{lemma}[theorem]{Lemma}
\newtheorem{corollary}[theorem]{Corollary}
\newtheorem{maintheorem}{Theorem}

\theoremstyle{definition}
\newtheorem{definition}[theorem]{Definition}
\newtheorem{example}[theorem]{Example}
\newtheorem{remark}[theorem]{Remark}

\numberwithin{equation}{section}

\newcommand{\C}{\mathbb C}
\newcommand{\Z}{\mathbb Z}
\newcommand{\id}{\mathrm{id}}
\newcommand{\cL}{\mathcal L}
\newcommand{\cE}{\mathcal E}
\newcommand{\cH}{\mathcal H}
\newcommand{\cM}{\mathcal M}
\newcommand{\YBCC}{\mathrm{YB}_{\mathrm{CC}}}
\newcommand{\EndCC}{\mathrm{End}_{\mathrm{CC}}}
\newcommand{\Diag}{\mathrm{Diag}}
\newcommand{\bsum}{\boxplus}
\DeclareMathOperator{\End}{End}
\DeclareMathOperator{\GL}{GL}
\DeclareMathOperator{\tr}{tr}
\DeclareMathOperator{\lk}{lk}
\DeclareMathOperator{\spn}{span}

\begin{document}

\title{Link invariants from charge-conserving Yang--Baxter operators}

\author{C\'esar Galindo}
\address{Departamento de Matem\'aticas, Universidad de los Andes, Bogot\'a, Colombia}
\email{cn.galindo1116@uniandes.edu.co}

\author{Laura Rodr\'iguez Sierra}
\address{Departamento de Matem\'aticas, Universidad de los Andes, Bogot\'a, Colombia}
\email{l.rodriguezs2@uniandes.edu.co}


\begin{abstract}
Turaev's enhanced Yang--Baxter operators produce invariants of oriented links from solutions of
the Yang--Baxter equation. If the enhancing endomorphism is diagonal with generic eigenvalues,
Turaev's compatibility condition forces the Yang--Baxter operator $R$ to be charge-conserving,
that is, to send each $e_i\otimes e_j$ into the span of $e_i\otimes e_j$ and $e_j\otimes e_i$.
Charge-conserving Yang--Baxter operators were classified by Martin and Rowell. Building on their
classification, we determine all enhancements of charge-conserving Yang--Baxter operators with
diagonal enhancing endomorphism, and we identify the resulting link invariants. The link
invariant associated with such an enhancement is a sum over the colourings of the components of
the link by the connected components of $R$, in which each monochromatic sublink contributes a
specialization of the HOMFLY--PT polynomial and different colours interact only through linking
numbers. Every such state sum arises in this way.
\end{abstract}

\maketitle

\section{Introduction}
\label{sec:introduction}

Let $V$ be a finite-dimensional complex vector space. A \emph{Yang--Baxter operator} on $V$
is an invertible $R\in\End(V^{\otimes2})$ such that
$(R\otimes\id)(\id\otimes R)(R\otimes\id)=(\id\otimes R)(R\otimes\id)(\id\otimes R)$.
Such an operator gives representations $\rho_{R,r}$ of the braid groups $B_r$ on
$V^{\otimes r}$, in which the generator $\sigma_k$ acts by $R$ on the tensor factors $k$ and
$k+1$. By the theorems of Alexander and Markov, an invariant of oriented links is
the same as a family of functions on the braid groups that is invariant under conjugation and
stabilization \cite{Ale,Mar}, \cite[Chapter~2]{KT}. Turaev~\cite{Tur} showed how to obtain such
functions from $R$. An \emph{enhancement} of $R$ consists of $\mu\in\End(V)$ and
$\alpha,\beta\in\C^\times$ satisfying a compatibility condition and two partial-trace
identities (see Section~\ref{subsec:turaev} for the precise definition), and then
\[
 \xi\longmapsto\alpha^{-w(\xi)}\beta^{-r}\tr\bigl(\rho_{R,r}(\xi)\,\mu^{\otimes r}\bigr)
 \qquad(\xi\in B_r),
\]
where $w(\xi)$ is the exponent sum, is an invariant of the closure of $\xi$. The
specializations of the HOMFLY--PT and Kauffman polynomials attached to the vector
representations of the quantum groups of classical type arise in this way
\cite[Section~4]{Tur}. Turaev remarked that ``the resources of this scheme are far from being
exhausted'' \cite[p.~527]{Tur}.

This paper grew out of an attempt to classify enhanced Yang--Baxter operators in a fixed
dimension. In dimension two all solutions of the Yang--Baxter equation are
known~\cite{Hie92,Hie93}, and Aizawa, Harada, Kawaguchi and Otsuki determined their
enhancements and invariants case by case~\cite{AHKO}; in higher dimensions no list of solutions
is available. One may instead organize the problem according to the enhancing endomorphism
$\mu$. A change of basis of $V$ changes neither the enhancement conditions nor the
invariant (Lemma~\ref{lem:conjugation}), and therefore we may assume that $\mu$ is in Jordan normal form.
Suppose that $\mu$ is diagonal, $\mu e_i=\mu_ie_i$ for a basis $(e_i)_{i\in I}$, and define the coefficients $R_{ij}^{kl}$ of $R$ by $R(e_i\otimes e_j)=\sum_{k,l}R_{ij}^{kl}\,e_k\otimes e_l$. Then the compatibility condition
$(\mu\otimes\mu)R=R(\mu\otimes\mu)$ reads
\[
 (\mu_k\mu_l-\mu_i\mu_j)\,R_{ij}^{kl}=0\qquad(i,j,k,l\in I)
\]
\cite[Theorem~2.3.1]{Tur}. If $\mu$ is \emph{generic}, in the sense that
$\mu_i\mu_j\neq\mu_k\mu_l$ whenever $\{i,j\}\neq\{k,l\}$ as multisets, this forces
$R_{ij}^{kl}=0$ unless $\{k,l\}=\{i,j\}$, that is,
\begin{equation}\label{eq:cc-intro}
 R(e_i\otimes e_j)\in\spn\{e_i\otimes e_j,\ e_j\otimes e_i\}\qquad(i,j\in I).
\end{equation}
Operators with this property are called \emph{charge-conserving}. Conversely, a diagonal
$\mu$ is compatible with every charge-conserving $R$, since $\mu\otimes\mu$ acts by the scalar
$\mu_i\mu_j$ on the span of $e_i\otimes e_j$ and $e_j\otimes e_i$. Thus, for diagonalizable $\mu$
with generic eigenvalues, the classification problem becomes the determination of the
\emph{diagonal enhancements} of charge-conserving Yang--Baxter operators, defined precisely
below.

Unlike the original problem, this one can be solved in all dimensions at once, because
charge-conserving Yang--Baxter operators were classified by Martin and Rowell~\cite{MR}, under
the name of spin-chain braid representations. We describe the outcome of their classification
in the next subsection. The study of the link invariants that Turaev's construction attaches to
these solutions is also among the directions for future work listed in \cite[Section~7]{MR}.

\subsection*{The data of a charge-conserving operator}
We now describe the combinatorial data that enter our results; the precise
statements, with references to~\cite{MR}, are in Section~\ref{sec:cc}. Let $R$ be a
charge-conserving Yang--Baxter operator with respect to a basis $(e_i)_{i\in I}$ of $V$.
By~\eqref{eq:cc-intro}, $R$ multiplies each $e_i\otimes e_i$ by a scalar, and for
$i\neq j$ it preserves the plane spanned by $e_i\otimes e_j$ and $e_j\otimes e_i$, on which it
acts by a $2\times2$ matrix $R_{ij}$, which we call the \emph{block} of the pair $\{i,j\}$. Martin and Rowell
show that these blocks are organized as follows (Propositions~\ref{prop:two-index}
and~\ref{prop:MR-structure}).
\begin{enumerate}[label=\textup{(\alph*)},leftmargin=*]
\item \emph{Connected components.} Call $\{i,j\}$ a \emph{flip edge} if both diagonal entries of
$R_{ij}$ vanish, that is, if $R$ maps $e_i\otimes e_j$ to a multiple of $e_j\otimes e_i$ and vice
versa. The relation ``$i=j$ or $\{i,j\}$ is not a flip edge'' is an equivalence relation on $I$.
We call its classes $\Lambda_1,\dots,\Lambda_m$ the \emph{connected components} of $R$, and we say that $R$ is \emph{connected} if $m=1$. Let $V_{\Lambda}$ denote the span of the $e_i$ with $i\in\Lambda$. Then
$R$ preserves $V_{\Lambda_s}\otimes V_{\Lambda_s}$, and its restriction $R^{(s)}:=R|_{V_{\Lambda_s}\otimes V_{\Lambda_s}}$ is a connected charge-conserving Yang--Baxter operator, which we call the \emph{connected restriction} of $R$ to $\Lambda_s$. Finally, if $i\in\Lambda_s$ and $j\in\Lambda_t$ with $s\neq t$, then $R_{ij}^2=\kappa_{st}\,\id$ for a scalar $\kappa_{st}\neq0$ that depends only on $s$ and $t$.
\item \emph{Scalar classes.} Let $R$ be connected. The relation ``$i=j$ or $R_{ij}$ is a scalar
matrix'' is again an equivalence relation. We call its classes the \emph{scalar classes} of $R$. On each class $C$ the operator $R$ multiplies every $e_i\otimes e_k$ with $i,k\in C$ by the same scalar $v(C)$, which we call the \emph{class scalar} of $C$. If $i\in C$ and $j\in D$ lie in different classes, then
$R$ maps exactly one of the two vectors $e_i\otimes e_j$ and $e_j\otimes e_i$ to a multiple of
the other; we write $C\prec D$ if it maps $e_j\otimes e_i$ to a multiple of $e_i\otimes e_j$.
This does not depend on the choice of $i\in C$ and $j\in D$, and it is a total order. We label
the classes $C_1\prec\dots\prec C_\ell$ and put $v_j:=v(C_j)$.
\item \emph{Spectral parameters.} Let $R$ be connected. If $\ell\ge2$, all blocks $R_{ij}$ with
$i$ and $j$ in different scalar classes have the same two eigenvalues $a,b$, possibly equal, with
$a+b\neq0$, and every class scalar is equal to $a$ or to $b$. If $\ell=1$, then $R=v_1\,\id$, and
we put $a:=v_1$ and $b:=v_1$. In both cases we call $a,b$ the \emph{spectral parameters} of $R$.
\end{enumerate}
For example, the standard
Hecke operator on $\C^N$, $N\ge2$ (Example~\ref{ex:hecke}), is connected, with scalar classes
$\{1\}\prec\dots\prec\{N\}$, spectral parameters $Q$ and $-Q^{-1}$, and all class scalars equal
to $Q$; and a scalar multiple $\lambda\tau$ of the flip $\tau(x\otimes y)=y\otimes x$ has the
singletons as connected components, with $\kappa_{st}=\lambda^2$.

\subsection*{Main results}
Absorbing Turaev's scalar $\beta$ into $\mu$ (Lemma~\ref{lem:turaev-relation}), a
\emph{diagonal enhancement} of a charge-conserving Yang--Baxter operator $R$ is a triple
$(R,\mu,\alpha)$, where $\alpha\in\C^\times$ and $\mu\in\GL(V)$ is diagonal in the basis
$(e_i)$, say $\mu e_i=\mu_ie_i$, such that
\[
 \tr_2\bigl(R(\id\otimes\mu)\bigr)=\alpha\,\id,\qquad
 \tr_2\bigl(R^{-1}(\id\otimes\mu)\bigr)=\alpha^{-1}\id;
\]
here $\tr_2\colon\End(V\otimes V)\to\End(V)$ is the partial trace over the second factor, and
Turaev's compatibility condition $(\mu\otimes\mu)R=R(\mu\otimes\mu)$ holds automatically. By
Turaev's theorem (Section~\ref{subsec:turaev}),
$T_{R,\mu,\alpha}(\widehat\xi)=\alpha^{-w(\xi)}\tr\bigl(\rho_{R,r}(\xi)\mu^{\otimes r}\bigr)$,
for $\xi\in B_r$, is then an invariant of the closure $\widehat\xi$ of $\xi$.

The following numbers, attached to a connected operator $R$, govern its enhancements. Fix an
ordering $(a,b)$ of the spectral parameters. For $0\le j\le\ell$ let $h_j$ be the number of
classes among $C_1,\dots,C_j$ with class scalar $a$ minus the number of those with class scalar
$b$ if $a\neq b$, and put $h_j:=j$ if $a=b$. We define the \emph{signed effective rank} of $R$ to be $n:=h_\ell$, and we set \[ \cE(R):=-ab\Bigl(-\frac ab\Bigr)^{n}.
\]
Exchanging $a$ and $b$ replaces every $h_j$ by $-h_j$, and it does not change $\cE(R)$
(Lemma~\ref{lem:v-identities}).

\begin{maintheorem}[Existence and parametrization; Theorem~\ref{thm:global-existence}]
\label{mainthm:A}
A charge-conserving Yang--Baxter operator $R$ with connected restrictions
$R^{(1)},\dots,R^{(m)}$ admits a diagonal enhancement if and only if
$\cE(R^{(1)})=\dots=\cE(R^{(m)})$. In that case, a triple $(R,\mu,\alpha)$ with
$\alpha\in\C^\times$ and $\mu\in\GL(V)$ diagonal is a diagonal enhancement if and only if
$\alpha^2$ is this common value and, for every scalar class $C_j$ of every $R^{(s)}$,
\begin{equation}\label{eq:class-traces-intro}
 \sum_{k\in C_j}\mu_k=\frac{v_j}{\alpha}\Bigl(-\frac ab\Bigr)^{h_{j-1}},
\end{equation}
where $a$, $b$, $v_j$ and $h_{j-1}$ are computed in $R^{(s)}$.
\end{maintheorem}

Thus, once $\alpha$ is fixed, only the sums of the eigenvalues of $\mu$ over the scalar classes
are prescribed; in particular $\mu$ is determined by $\alpha$ when all scalar classes are
singletons.

For $Q\in\C^\times$ with $Q^2\neq1$ and $k\in\Z$ let $[k]_Q:=(Q^k-Q^{-k})/(Q-Q^{-1})$, and let
$P_L(A,z)$ be the HOMFLY--PT polynomial of an oriented link $L$, normalized by
\begin{equation}\label{eq:homfly-intro}
 P_{\bigcirc}=1,\qquad A\,P_{L_+}-A^{-1}P_{L_-}=z\,P_{L_0},
\end{equation}
where $\bigcirc$ is the unknot and $L_+,L_-,L_0$ are links with diagrams that coincide outside a
small disc, inside which they are a positive crossing, a negative crossing and its oriented
smoothing (Section~\ref{subsec:homfly}).

\begin{maintheorem}[Connected invariants; Theorem~\ref{thm:connected-invariant}]
\label{mainthm:B}
Let $(R,\mu,\alpha)$ be a diagonal enhancement of a connected $R$ with spectral parameters
$a,b$ and signed effective rank $n$. Choose $d\in\C^\times$ with $d^2=-ab$ and put $Q:=a/d$.
Then $Q^2\neq1$, $\varepsilon:=\alpha/(dQ^n)\in\{\pm1\}$, and for every oriented link $L$
\[
 T_{R,\mu,\alpha}(L)=\varepsilon\,[n]_Q\;P_L\bigl(\varepsilon Q^n,\;Q-Q^{-1}\bigr).
\]
\end{maintheorem}

The right-hand side does not depend on the choice of $d$ or of the ordering $(a,b)$
(Remark~\ref{rem:choices}). Thus the invariant depends only on $Q$, $n$ and the sign $\varepsilon$. The dimension of $V$, the sizes and the order of the scalar classes, and the individual eigenvalues of $\mu$ are invisible to it. For the standard Hecke operator one has $n=\dim V$, and
Theorem~\ref{mainthm:B} recovers Turaev's construction of the $\mathfrak{sl}_N$
specializations of the HOMFLY--PT polynomial \cite[Section~4.2]{Tur}.

Now let $R$ be arbitrary, with connected components $\Lambda_1,\dots,\Lambda_m$. If
$(R,\mu,\alpha)$ is a diagonal enhancement, then $(R^{(s)},\mu^{(s)},\alpha)$ is also a diagonal enhancement for every
$s$, where $\mu^{(s)}$ is the restriction of $\mu$ to $V_{\Lambda_s}$
(Proposition~\ref{prop:decoupling}). A \emph{colouring} of a link $L$ is a map
$\chi\colon\pi_0(L)\to\{1,\dots,m\}$ from its set of components; $L^\chi_s$ denotes the sublink
of components of colour $s$, and $\lk$ the linking number, with $\lk(J,\varnothing):=0$.

\begin{maintheorem}[State sum; Theorem~\ref{thm:state-sum}]
\label{mainthm:C}
Let $(R,\mu,\alpha)$ be a diagonal enhancement, and let $(Q_s,n_s,\varepsilon_s)$ be the data of
Theorem~\ref{mainthm:B} for $(R^{(s)},\mu^{(s)},\alpha)$. Put $A_s:=\varepsilon_sQ_s^{n_s}$, $z_s:=Q_s-Q_s^{-1}$, $\delta_s:=\varepsilon_s[n_s]_{Q_s}$ and $\omega_{st}:=\kappa_{st}/\alpha^2$.
Then for every oriented link $L$
\begin{equation}\label{eq:state-sum-intro}
 T_{R,\mu,\alpha}(L)=\sum_{\chi\colon\pi_0(L)\to\{1,\dots,m\}}\
 \prod_{s<t}\omega_{st}^{\lk(L^\chi_s,L^\chi_t)}
 \prod_{s\,:\,L^\chi_s\neq\varnothing}\delta_s\,P_{L^\chi_s}(A_s,z_s).
\end{equation}
\end{maintheorem}

\begin{maintheorem}[Realization; Theorem~\ref{thm:realization}]
\label{mainthm:D}
Let $m\ge1$, and for $1\le s\le m$ let $Q_s\in\C^\times$ with $Q_s^2\neq1$, $n_s\in\Z$ and
$\varepsilon_s\in\{\pm1\}$; for $s\neq t$ let $\omega_{st}\in\C^\times$ with $\omega_{st}=\omega_{ts}$. Then the
right-hand side of~\eqref{eq:state-sum-intro} is the invariant $T_{R,\mu,1}$ of a diagonal
enhancement of some charge-conserving Yang--Baxter operator. Consequently, the link invariants
of diagonally enhanced charge-conserving Yang--Baxter operators are exactly these state sums.
\end{maintheorem}

In particular, on knots these invariants are linear combinations
$\sum_s\delta_sP_K(A_s,z_s)$ of HOMFLY--PT specializations, and on links they are determined
by the HOMFLY--PT polynomials of all sublinks together with the pairwise linking numbers
(Corollary~\ref{cor:state-sum-consequences}). By the discussion at the beginning of this
introduction, together with Lemmas~\ref{lem:conjugation} and~\ref{lem:turaev-relation},
Theorems~\ref{mainthm:A}--\ref{mainthm:D} apply to every enhanced Yang--Baxter operator whose
enhancing endomorphism is diagonalizable with nonzero generic eigenvalues.

The paper is organized as follows. Section~\ref{sec:prelim} collects the facts on braids, the HOMFLY--PT polynomial and
Turaev's construction that we need, and describes the combinatorics of coloured braids used in
Section~\ref{sec:invariants}. Section~\ref{sec:cc} recalls from~\cite{MR} what we need on
charge-conserving Yang--Baxter operators. Section~\ref{sec:enhancements} shows that the
X-transformations used in~\cite{MR} do not change the link invariants, and determines all
diagonal enhancements (Theorem~\ref{mainthm:A}). Section~\ref{sec:invariants} computes their invariants
and proves Theorems~\ref{mainthm:B}--\ref{mainthm:D}. Section~\ref{sec:examples} compares these
results with the literature, namely the solutions in dimension two and their dressings
studied in~\cite{AHKO} and the invariants obtained from the Yokonuma--Hecke algebras, and ends
with some remarks.

\subsection*{Conventions}
All vector spaces are finite-dimensional over $\C$. Links are oriented, nonempty and
considered up to ambient isotopy in $\mathbb R^3$; $\pi_0(L)$ is the set of components of $L$
and $|L|=|\pi_0(L)|$ is their number. We write $\bigcirc$ for the unknot, $\bigcirc_r$ for the
$r$-component unlink and $L\sqcup L'$ for the split union. Braid diagrams are read from bottom
to top, and the braid generator $\sigma_k$ is a positive crossing
(Section~\ref{subsec:braids}).

\section{Braids, the HOMFLY--PT polynomial and Turaev's construction}
\label{sec:prelim}

This section fixes our conventions and collects the facts that we need. Most of it is
standard. The exception is Section~\ref{subsec:colourings}, which explains how to follow
colours along the strands of a braid. These combinatorial facts are the main ingredient in the
proof of the state sum formula in Section~\ref{sec:invariants}.

\subsection{Braids and closures}
\label{subsec:braids}
For $r\ge1$, the braid group $B_r$ has generators $\sigma_1,\dots,\sigma_{r-1}$ and relations
$\sigma_k\sigma_{k+1}\sigma_k=\sigma_{k+1}\sigma_k\sigma_{k+1}$ and
$\sigma_k\sigma_j=\sigma_j\sigma_k$ for $|k-j|\ge2$. In particular, $B_1$ is the trivial group. The \emph{exponent sum}
$w\colon B_r\to\Z$ and the \emph{underlying permutation} $\xi\mapsto\pi_\xi\in\mathfrak S_r$
are the homomorphisms with $w(\sigma_k)=1$ and $\pi_{\sigma_k}=(k\ k{+}1)$, the
transposition of $k$ and $k+1$.

We fix the following conventions for diagrams, because in Section~\ref{sec:invariants} we
shall follow tensor factors along the strands of a braid. Let
$u=\sigma_{k_1}^{\epsilon_1}\cdots\sigma_{k_q}^{\epsilon_q}$, with $\epsilon_p=\pm1$, be a word for $\xi\in B_r$.
Its \emph{braid diagram} $D_u$ consists of $r$ strands, oriented upwards, which join $r$ bottom
endpoints to $r$ top endpoints; both sets of endpoints are numbered $1,\dots,r$ from left to
right. The diagram has $q$ crossings, one for each letter, stacked in such a way that the last letter
$\sigma_{k_q}^{\epsilon_q}$ gives the lowest crossing and the first letter $\sigma_{k_1}^{\epsilon_1}$ the
highest one. The letter $\sigma_k^{\pm1}$ is a crossing of sign $\pm1$ between the strands at
positions $k$ and $k+1$; in particular, $\sigma_k$ is a positive crossing. Reading the word
from right to left thus corresponds to reading the diagram from bottom to top, which matches
the convention $\rho_{R,r}(\xi)=R_{k_1}^{\epsilon_1}\cdots R_{k_q}^{\epsilon_q}$ of
Section~\ref{subsec:turaev}, in which vectors enter at the bottom.

For $0\le p\le q$, \emph{level $p$} is a horizontal line that separates the $p$ highest
crossings from the others. Thus level $q$ is the bottom, level $0$ is the top, and the $p$-th
crossing lies between the levels $p$ and $p-1$. Passing the $p$-th crossing upwards exchanges
the positions $k_p$ and $k_p+1$ and leaves the other positions unchanged. Following a strand
from the bottom to the top therefore applies the transpositions
$\pi_{\sigma_{k_q}},\pi_{\sigma_{k_{q-1}}},\dots,\pi_{\sigma_{k_1}}$ to its position. Hence the strand that starts at the bottom position $i$ ends at the top position
$\pi_{\sigma_{k_1}}\cdots\pi_{\sigma_{k_q}}(i)=\pi_\xi(i)$.

The \emph{closure} $\widehat D_u$ is obtained by joining, for every $i$, the top endpoint and
the bottom endpoint numbered $i$ by an arc running to the right of $D_u$, the arcs being
disjoint. It is a diagram of the closed braid $\widehat\xi$, and its crossings are those of
$D_u$. Starting at the bottom position $i$ and travelling along $\widehat D_u$, we reach the top
position $\pi_\xi(i)$ and return to the bottom position $\pi_\xi(i)$; hence the components of
$\widehat\xi$ correspond to the cycles of $\pi_\xi$, the component through the bottom position
$i$ passing through the bottom positions $i,\pi_\xi(i),\pi_\xi^2(i),\dots$
\cite[Exercise~2.2.1]{KT}. The permutation $\pi_\xi$ is a product of $q$ transpositions and
$q\equiv w(\xi)\pmod2$, hence its sign is $(-1)^{w(\xi)}$. A permutation of $r$ letters with $k$
cycles has sign $(-1)^{r-k}$, and $\pi_\xi$ has $|\widehat\xi|$ cycles. Therefore
\begin{equation}\label{eq:parity}
 r-w(\xi)\equiv|\widehat\xi|\pmod 2 .
\end{equation}

Every oriented link is isotopic to a closed braid (Alexander's theorem~\cite{Ale}; see
\cite[Theorem~2.3]{KT}). Two braids $\xi\in B_r$ and $\xi'\in B_{r'}$ have isotopic closures if
and only if they are related by a finite sequence of conjugations in a fixed braid group and
stabilizations $\xi\leftrightarrow\iota(\xi)\sigma_r^{\pm1}$, where $\iota\colon B_r\to B_{r+1}$
is the natural inclusion (Markov's theorem~\cite{Mar}; see \cite[Theorem~2.8]{KT}
and~\cite{Mor}). In \cite{KT} stabilization is written $\xi\mapsto\sigma_r^{\pm1}\iota(\xi)$,
which is conjugate to $\iota(\xi)\sigma_r^{\pm1}$. Consequently, a family of maps
$f_r\colon B_r\to\C$ that is invariant under conjugation and stabilization (a \emph{Markov
function} \cite[Definition~2.10]{KT}) defines an invariant of oriented links by
$L=\widehat\xi\mapsto f_r(\xi)$ \cite[Corollary~2.9]{KT}.

For disjoint sublinks $J,K$ of a link $L$, $\lk(J,K)$ denotes half the sum of the signs of
the crossings between $J$ and $K$ in a diagram of $L$. It is the sum of the linking numbers of
the pairs formed by a component of $J$ and a component of $K$, and we set $\lk(J,\varnothing):=0$.

\subsection{Colourings of closed braids}
\label{subsec:colourings}
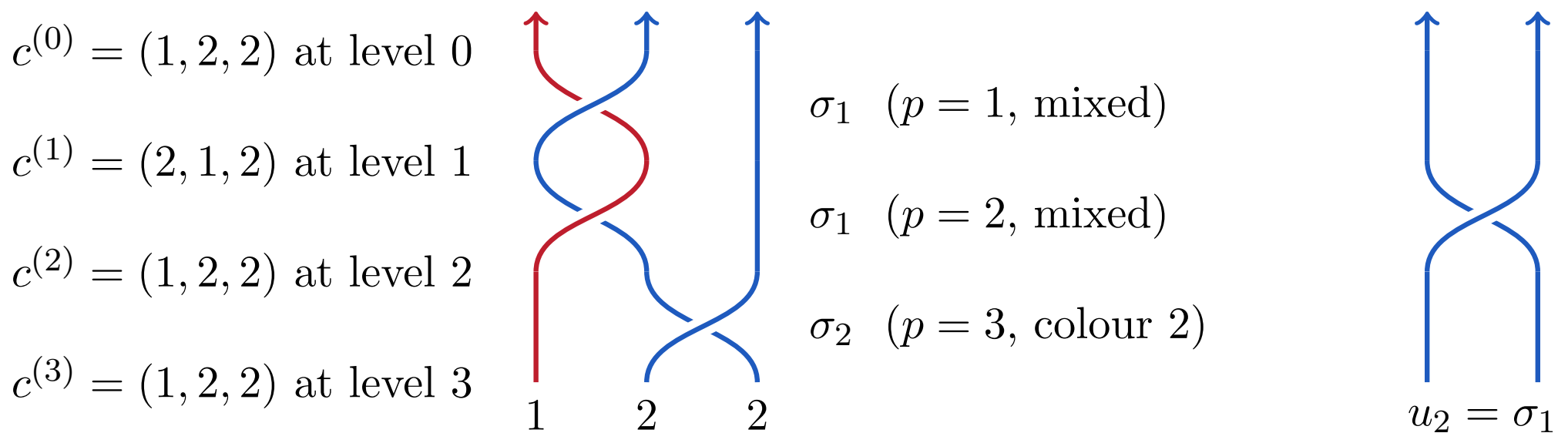
\begin{figure}[ht]
\centering
\begin{tikzpicture}[x=0.95cm,y=0.95cm,line width=1.2pt]
 \definecolor{cA}{RGB}{190,30,45}
 \definecolor{cB}{RGB}{30,90,190}
 % over strand from (x,y) to (x+1,y+1), under strand from (x+1,y) to (x,y+1)
 \newcommand{\poscross}[4]{%
  \draw[#4] (#1+1,#2) .. controls (#1+1,#2+0.5) and (#1,#2+0.5) .. (#1,#2+1);
  \draw[white,line width=5pt] (#1,#2) .. controls (#1,#2+0.5) and (#1+1,#2+0.5) .. (#1+1,#2+1);
  \draw[#3] (#1,#2) .. controls (#1,#2+0.5) and (#1+1,#2+0.5) .. (#1+1,#2+1);}
 % p=3: sigma_2 between levels 3 (y=0) and 2 (y=1)
 \draw[cA] (1,0) -- (1,1);
 \poscross{2}{0}{cB}{cB}
 % p=2: sigma_1 between y=1 and y=2; colours at level 2: (1,2,2)
 \poscross{1}{1}{cA}{cB}
 \draw[cB] (3,1) -- (3,2);
 % p=1: sigma_1 between y=2 and y=3; colours at level 1: (2,1,2)
 \poscross{1}{2}{cB}{cA}
 \draw[cB] (3,2) -- (3,3);
 % arrows indicating orientation
 \foreach \x/\c in {1/cA,2/cB,3/cB} {\draw[\c,->] (\x,3) -- (\x,3.35);}
 % colours at the bottom endpoints
 \node[below] at (1,0) {\small $1$};
 \node[below] at (2,0) {\small $2$};
 \node[below] at (3,0) {\small $2$};
 % colour words at the levels
 \node[left] at (0.6,0) {\small $c^{(3)}=(1,2,2)$ at level $3$};
 \node[left] at (0.6,1) {\small $c^{(2)}=(1,2,2)$ at level $2$};
 \node[left] at (0.6,2) {\small $c^{(1)}=(2,1,2)$ at level $1$};
 \node[left] at (0.6,3) {\small $c^{(0)}=(1,2,2)$ at level $0$};
 % letters
 \node[right] at (3.3,0.5) {\small $\sigma_2$\ \ ($p=3$, colour $2$)};
 \node[right] at (3.3,1.5) {\small $\sigma_1$\ \ ($p=2$, mixed)};
 \node[right] at (3.3,2.5) {\small $\sigma_1$\ \ ($p=1$, mixed)};
 % the deleted braid u_2 = sigma_1
 \begin{scope}[xshift=8.6cm]
  \draw[cB] (0,0) -- (0,1); \draw[cB] (1,0) -- (1,1);
  \poscross{0}{1}{cB}{cB}
  \draw[cB] (0,2) -- (0,3); \draw[cB] (1,2) -- (1,3);
  \foreach \x in {0,1} {\draw[cB,->] (\x,3) -- (\x,3.35);}
  \node[below] at (0.5,-0.05) {\small $u_2=\sigma_1$};
 \end{scope}
\end{tikzpicture}
\caption{The diagram $D_u$ of $u=\sigma_1\sigma_1\sigma_2\in B_3$, read from bottom to top,
with the compatible colour word $c=(1,2,2)$. Colour $1$ is red and colour $2$ is blue, and the
colour of each strand is also written below its bottom endpoint. The colour words at the
levels are shown on the left. Deleting the strand of colour $1$ leaves the diagram of
$u_2=\sigma_1\in B_2$ (right).}
\label{fig:coloured-braid}
\end{figure}

In Section~\ref{sec:invariants} the basis vectors of $V$ will carry colours, and a braid will
move each tensor factor along a strand without changing its colour. We now set up the
corresponding combinatorics; Figure~\ref{fig:coloured-braid} illustrates all the notions
introduced here.

Let $u=\sigma_{k_1}^{\epsilon_1}\cdots\sigma_{k_q}^{\epsilon_q}$ be a word for $\xi\in B_r$, and let $m\ge1$.
The symmetric group $\mathfrak S_r$ acts on \emph{colour words} $c=(c_1,\dots,c_r)\in
\{1,\dots,m\}^r$ by $(\pi\cdot c)_i=c_{\pi^{-1}(i)}$; thus $\pi\cdot c$ is obtained from $c$ by
moving the entry in position $i$ to position $\pi(i)$. We say that $c$ is \emph{compatible}
with $\xi$ if $\pi_\xi\cdot c=c$, that is, if $c$ is constant on the cycles of $\pi_\xi$. Since
the cycles of $\pi_\xi$ correspond to the components of $\widehat\xi$, a compatible colour word
is the same as a colouring of the components of $\widehat\xi$, in which the entry $c_i$ is the
colour of the component through the bottom position $i$.

Fix a compatible colour word $c$ and colour every strand of $D_u$ by the colour of its bottom
endpoint. For $0\le p\le q$ let $c^{(p)}$ be the word formed by the colours of the strands at
level $p$, read from left to right. Since the $p$-th crossing exchanges the positions $k_p$
and $k_p+1$,
\[
 c^{(q)}=c,\qquad c^{(p-1)}=\pi_{\sigma_{k_p}}\cdot c^{(p)}\quad(1\le p\le q),\qquad
 c^{(0)}=\pi_\xi\cdot c=c .
\]
The last equality says that the colours at the top agree with the colours at the bottom, as
they must, since the closure joins the top endpoint $i$ to the bottom endpoint $i$.

The $p$-th crossing joins the strands at the positions $k_p$ and $k_p+1$ of level $p$. It is
\emph{monochromatic of colour $s$} if both strands have colour $s$, that is, if
$c^{(p)}_{k_p}=c^{(p)}_{k_p+1}=s$, and \emph{mixed} otherwise. For a colour $s$ put $S_s:=\{i:c_i=s\}$. At every level exactly $|S_s|$ strands have colour $s$, and we number them
from left to right. If the $p$-th crossing is monochromatic of colour $s$, its
two strands are adjacent in this numbering; the left one has number \[ j_p:=\#\{i\le k_p:\ c^{(p)}_i=s\},
\]
and the right one has number $j_p+1$. Let $u_s$ be the word obtained from $u$ by deleting the
letters of all crossings that are not monochromatic of colour $s$, and replacing every
remaining letter $\sigma_{k_p}^{\epsilon_p}$ by $\sigma_{j_p}^{\epsilon_p}$. Informally, $u_s$ is the braid
formed by the strands of colour $s$ alone.

\begin{example}\label{ex:coloured-braid}
Let $u=\sigma_1\sigma_1\sigma_2\in B_3$ (Figure~\ref{fig:coloured-braid}). Here
$k_1=k_2=1$, $k_3=2$ and $\pi_\xi$ is the transposition $(2\ 3)$. Hence the closure has two components, namely
one through the bottom position $1$ and one through the bottom positions $2$ and $3$. The colour
word $c=(1,2,2)$ is compatible with $\xi$, and $c^{(3)}=c^{(2)}=(1,2,2)$,
$c^{(1)}=(2,1,2)$, $c^{(0)}=(1,2,2)$. The lowest crossing ($p=3$) is monochromatic of colour
$2$, and the other two crossings are mixed. Hence $u_1$ is the empty word in $B_1$, and
$u_2=\sigma_1\in B_2$, because $j_3=1$. Both components are unknots, and the two mixed
crossings are positive; indeed $\widehat\xi$ is the positive Hopf link, the closure of
$\sigma_1^2\in B_2$, obtained from $\xi=\iota(\sigma_1^2)\sigma_2$ by one stabilization.
\end{example}

\begin{lemma}\label{lem:deletion}
Let $u$ be a word for $\xi\in B_r$ and $c$ a colour word compatible with $\xi$, as above. For
a colour $s$, let $L_s\subseteq\widehat\xi$ be the sublink formed by the components of
colour $s$, that is, by the components through the bottom positions in $S_s$
\textup{(}$L_s=\varnothing$ if $S_s=\varnothing$\textup{)}, and for $s\neq t$ let $w_{st}$ be
the sum of the signs of the mixed crossings between strands of colours $s$ and $t$. Then the
following hold.
\begin{enumerate}[label=\textup{(\roman*)}]
\item If $S_s\neq\varnothing$, then $u_s$ is a word for a braid in $B_{|S_s|}$ whose closure is
isotopic to $L_s$.
\item $w_{st}=2\lk(L_s,L_t)$.
\item $w(\xi)=\sum_{s\,:\,S_s\neq\varnothing}w(u_s)+\sum_{s<t}w_{st}$.
\end{enumerate}
\end{lemma}

In Example~\ref{ex:coloured-braid}, (i) says that both components are unknots, and (ii)
and (iii) read $w_{12}=2=2\lk(L_1,L_2)$ and $w(\xi)=3=0+1+2$.

\begin{proof}
(i) The strands of colour $s$ start at the bottom positions in $S_s$ and end at the top
positions $\pi_\xi(S_s)=S_s$. Delete from $D_u$ all strands whose colour is not $s$. This removes
exactly the crossings that are not monochromatic of colour $s$, and what remains is a braid
diagram with $|S_s|$ strands, still oriented upwards. Renumber the remaining positions at each
level from left to right. The $p$-th crossing, which joined the positions $k_p$ and $k_p+1$,
now joins the positions $j_p$ and $j_p+1$, and its sign is unchanged. Hence the new diagram is
$D_{u_s}$. Its top endpoints and its bottom endpoints both occupy the original positions $S_s$,
renumbered in increasing order, and therefore closing $D_{u_s}$ amounts to keeping the arcs of
$\widehat D_u$ that join the endpoints in $S_s$. Therefore $\widehat{D_{u_s}}$ is obtained from
$\widehat D_u$ by deleting the components of colours other than $s$, and it is a diagram of
$L_s$.

(ii) The crossings of $\widehat D_u$ between $L_s$ and $L_t$ are exactly the mixed crossings
between strands of colours $s$ and $t$, and $\lk(L_s,L_t)$ is half the sum of their signs.

(iii) The exponent sum $w(\xi)$ is the sum of the signs of all crossings of $D_u$, and $w(u_s)$
is the sum of the signs of the monochromatic crossings of colour $s$. Sorting the crossings
according to the colours of their two strands gives the formula.
\end{proof}

\subsection{The HOMFLY--PT polynomial}
\label{subsec:homfly}
Three oriented links $L_+,L_-,L_0$ form a \emph{skein triple} if they have diagrams that
coincide outside a disc, inside which they are a positive crossing, a negative crossing and its
oriented smoothing. The HOMFLY--PT polynomial $P_L(A,z)\in\Z[A^{\pm1},z^{\pm1}]$ is the unique
invariant of oriented links such that $P_\bigcirc=1$ and $A\,P_{L_+}-A^{-1}P_{L_-}=z\,P_{L_0}$
for every skein triple~\cite{HOMFLY,PT}; this is the normalization~\eqref{eq:homfly-intro}. If $L_+$ and $L_-$ are diagrams of a
link $L$ with a positive and a negative kink added, then $L_0$ is a diagram of $L\sqcup\bigcirc$. Therefore $(A-A^{-1})P_L=zP_{L\sqcup\bigcirc}$, and hence
$P_{\bigcirc_r}=\bigl(\frac{A-A^{-1}}{z}\bigr)^{r-1}$. The substitution $(A,z)\mapsto(-A,-z)$
preserves~\eqref{eq:homfly-intro}, and uniqueness gives
\begin{equation}\label{eq:homfly-sign}
 P_L(-A,-z)=P_L(A,z).
\end{equation}
We evaluate $P_L$ at points $(A,z)\in(\C^\times)^2$. To identify an invariant with a
specialization of $P_L$ we use the following standard uniqueness statement; its proof is the
descending-diagram argument of~\cite{LM}, which we recall for completeness.

\begin{lemma}\label{lem:skein-uniqueness}
Let $A\in\C^\times$ and $z\in\C$. Let $\Psi,\Psi'$ be $\C$-valued invariants of oriented links such
that $A\,\Psi(L_+)-A^{-1}\Psi(L_-)=z\,\Psi(L_0)$ for every skein triple, and the same for $\Psi'$.
If $\Psi(\bigcirc_r)=\Psi'(\bigcirc_r)$ for all $r\ge1$, then $\Psi=\Psi'$.
\end{lemma}

\begin{proof}
The difference $F=\Psi-\Psi'$ satisfies the same skein relation and vanishes on unlinks. We show
that $F$ vanishes on every link. For a link diagram $D$ we write $F(D)$ for the value of $F$
on the link represented by $D$.

Order the components of $D$, choose a base point on each of them away from the crossings, and
traverse the components in this order, each one from its base point along its orientation.
Call a crossing \emph{bad} if it is first reached along its under-strand. A diagram without bad
crossings is called \emph{descending}. In such a diagram each component lies below the previous
ones and descends from its base point, and therefore the diagram represents an unlink \cite{LM}.

We argue by induction on the number $n(D)$ of crossings and, for fixed $n(D)$, on the number
$b(D)$ of bad crossings. If $b(D)=0$, then $D$ represents an unlink and $F(D)=0$; this covers
$n(D)=0$. If $b(D)>0$, let $x$ be a bad crossing, and let $D'$ and $D_0$ be obtained from $D$ by
switching and by smoothing $x$. With the same order and base points, $D'$ has the same
components as $D$, the same crossings, and one bad crossing less. Thus $n(D')=n(D)$ and
$b(D')=b(D)-1$; and $n(D_0)=n(D)-1$. The skein relation, applied to the triple formed by $D$,
$D'$ and $D_0$, gives
\[
 F(D)=A^{\mp2}F(D')\pm A^{\mp1}z\,F(D_0)
\]
according as $x$ is positive or negative. Both terms on the right vanish by induction.
\end{proof}

\subsection{Turaev's construction}
\label{subsec:turaev}
For $F\in\End(V^{\otimes2})$ and $1\le k<r$ put $F_k:=\id^{\otimes(k-1)}\otimes F\otimes\id^{\otimes(r-k-1)}\in\End(V^{\otimes r})$. If $R$ is a
Yang--Baxter operator, then $\sigma_k\mapsto R_k$ defines a representation
$\rho_{R,r}\colon B_r\to\GL(V^{\otimes r})$, and
$\rho_{R,r}(\xi)=R_{k_1}^{\epsilon_1}\cdots R_{k_q}^{\epsilon_q}$ for every word
$\sigma_{k_1}^{\epsilon_1}\cdots\sigma_{k_q}^{\epsilon_q}$ for $\xi$. The partial trace
$\tr_2\colon\End(V^{\otimes2})\to\End(V)$ is
$\tr_2(F)(x)=\sum_l(\id\otimes e^l)F(x\otimes e_l)$ for a basis $(e_l)$ of $V$ with dual
basis $(e^l)$; it does not depend on the basis.

An \emph{enhanced Yang--Baxter operator} \cite[Section~2.3]{Tur} is a quadruple
$S=(R,\mu,\alpha,\beta)$, where $R$ is a Yang--Baxter operator, $\mu\in\End(V)$ and
$\alpha,\beta\in\C^\times$, such that
\begin{enumerate}[label=\textup{(\roman*)}]
\item $(\mu\otimes\mu)R=R(\mu\otimes\mu)$, and
\item $\tr_2\bigl(R^{\pm1}(\mu\otimes\mu)\bigr)=\alpha^{\pm1}\beta\,\mu$.
\end{enumerate}
Since $\tr_2\bigl(F(\mu\otimes\mu)\bigr)=\tr_2\bigl(F(\id\otimes\mu)\bigr)\mu$, if $\mu$ is
invertible then (ii) is equivalent to $\tr_2\bigl(R^{\pm1}(\id\otimes\mu)\bigr)=
\alpha^{\pm1}\beta\,\id$ \cite[Section~2.3]{Tur}. Condition (i) implies that $\mu^{\otimes r}$
commutes with $\rho_{R,r}(\xi)$, which makes the trace below invariant under conjugation of
$\xi$; condition (ii) implies that a stabilization $\xi\mapsto\iota(\xi)\sigma_r^{\pm1}$
multiplies the trace by $\alpha^{\pm1}\beta$, which is compensated by the normalizing factor.
This is the idea of the proof of the following theorem.

\begin{theorem}[{\cite[Theorem~3.1.2 and Section~3.2]{Tur}}]\label{thm:turaev}
If $S=(R,\mu,\alpha,\beta)$ is an enhanced Yang--Baxter operator, then
$\xi\mapsto\alpha^{-w(\xi)}\beta^{-r}\tr\bigl(\rho_{R,r}(\xi)\mu^{\otimes r}\bigr)$ is a Markov
function, and hence defines an invariant $T_S$ of oriented links. Moreover
$T_S(\bigcirc)=\beta^{-1}\tr(\mu)$ and $T_S(L\sqcup L')=T_S(L)T_S(L')$.
\end{theorem}

\begin{theorem}[{\cite[Theorem~3.2.1]{Tur}}]\label{thm:turaev-skein}
Let $S=(R,\mu,\alpha,\beta)$ be an enhanced Yang--Baxter operator, and let $\gamma_i\in\C$,
$i\in\Z$, be scalars, only finitely many of them nonzero, with $\sum_i\gamma_iR^i=0$. If the
links $L^{(i)}$, for the $i$ with $\gamma_i\neq0$, have diagrams that coincide outside a disc,
inside which $L^{(i)}$ is the two-strand braid $\sigma_1^i$, then
$\sum_i\gamma_i\alpha^iT_S(L^{(i)})=0$. In particular, if $(R-a)(R-b)=0$ with $ab\neq0$,
then for every skein triple
\begin{equation}\label{eq:quadratic-skein}
 \alpha\,T_S(L_+)-(a+b)\,T_S(L_0)+ab\,\alpha^{-1}\,T_S(L_-)=0 .
\end{equation}
\end{theorem}

The second statement follows from the first applied to the relation
$R-(a+b)\id+abR^{-1}=0$, which is $(R-a)(R-b)=0$ multiplied by $R^{-1}$, with
$L^{(1)}=L_+$, $L^{(0)}=L_0$ and $L^{(-1)}=L_-$.

Finally, a change of basis of $V$ does not affect enhancements or invariants.

\begin{lemma}\label{lem:conjugation}
Let $S=(R,\mu,\alpha,\beta)$ be an enhanced Yang--Baxter operator on $V$ and $g\in\GL(V)$, and
put $R^g:=(g\otimes g)R(g\otimes g)^{-1}$. Then $S^g:=\bigl(R^g,\,g\mu g^{-1},\,\alpha,\,\beta\bigr)$ is an enhanced Yang--Baxter operator, and
$T_{S^g}=T_S$.
\end{lemma}

\begin{proof}
We have $(R^g)_k=g^{\otimes r}R_k(g^{\otimes r})^{-1}$, hence $R^g$ is a Yang--Baxter operator and
$\rho_{R^g,r}(\xi)=g^{\otimes r}\rho_{R,r}(\xi)(g^{\otimes r})^{-1}$. Condition (i) for $S^g$ is
condition (i) for $S$ conjugated by $g\otimes g$. For condition (ii) we use
$\tr_2\bigl((g\otimes g)F(g\otimes g)^{-1}\bigr)=g\,\tr_2(F)\,g^{-1}$, which follows from the
two identities $\tr_2\bigl((\id\otimes g)F(\id\otimes g)^{-1}\bigr)=\tr_2(F)$ and
$\tr_2\bigl((g\otimes\id)F(g\otimes\id)^{-1}\bigr)=g\,\tr_2(F)\,g^{-1}$. Finally,
$T_{S^g}=T_S$ by the cyclicity of the trace.
\end{proof}

\section{Charge-conserving Yang--Baxter operators}
\label{sec:cc}

This section collects what we need from the classification of Martin and Rowell~\cite{MR},
translated into our conventions; the results attributed to~\cite{MR} are proved there. Our
matrices act on column vectors, and therefore the block $A(i,j)$ of \cite[(16)]{MR} is the transpose of
the block~\eqref{eq:edge-block} below. Proposition~\ref{thm:canonical-decomposition} and
Corollary~\ref{cor:quadratic} are short consequences of these results. The charge-conserving
representations of the loop braid groups were classified in the same spirit by Martin, Rowell
and Torzewska~\cite{MRT}.

\subsection{Coefficients and blocks}
\label{subsec:blocks}
An \emph{indexed line decomposition} of $V$ is a family $\cL=(L_i)_{i\in I}$ of
one-dimensional subspaces, indexed by a finite nonempty set $I$, with
$V=\bigoplus_{i\in I}L_i$. For $\varnothing\neq J\subseteq I$ let $V_J:=\bigoplus_{i\in J}L_i$ and $\cL_J:=(L_i)_{i\in J}$. For $i\neq j$ put $W_{ij}:=(L_i\otimes L_j)\oplus(L_j\otimes L_i)$. Then
\begin{equation}\label{eq:charge-decomposition}
 V^{\otimes2}=\bigoplus_{i\in I}(L_i\otimes L_i)\oplus\bigoplus_{\{i,j\}\subseteq I,\ i\neq j}W_{ij}.
\end{equation}
An operator on $V^{\otimes2}$ is \emph{charge-conserving} if it preserves every summand
of~\eqref{eq:charge-decomposition}; these operators form an algebra
$\EndCC(V^{\otimes2};\cL)$. We write $\YBCC(V,\cL)$ for the set of invertible
charge-conserving operators satisfying the Yang--Baxter equation.

Choose $0\neq e_i\in L_i$ and define the coefficients of $F\in\End(V^{\otimes2})$ by
$F(e_i\otimes e_j)=\sum_{k,l}F_{ij}^{kl}\,e_k\otimes e_l$. Then $F$ is charge-conserving if and
only if $F_{ij}^{kl}=0$ unless $(k,l)\in\{(i,j),(j,i)\}$, which is
condition~\eqref{eq:cc-intro}. Replacing $e_i$ by $\lambda_ie_i$ multiplies $F_{ij}^{kl}$ by $\lambda_i\lambda_j/(\lambda_k\lambda_l)$, and therefore the coefficients of a charge-conserving operator do
not depend on the choice of the $e_i$. For $R\in\EndCC(V^{\otimes2};\cL)$ we call the scalars $v_i:=R_{ii}^{ii}$, $i\in I$, the \emph{vertex scalars} of $R$. For $i\neq j$ we call the restriction $R_{ij}:=R|_{W_{ij}}$ the \emph{block} of the pair $\{i,j\}$; it has the matrix
\begin{equation}\label{eq:edge-block}
 R_{ij}=\begin{pmatrix}R_{ij}^{ij}&R_{ji}^{ij}\\[2pt] R_{ij}^{ji}&R_{ji}^{ji}\end{pmatrix}
\end{equation}
in the ordered basis $(e_i\otimes e_j,\,e_j\otimes e_i)$; reversing the order of $i,j$
conjugates it by $\bigl(\begin{smallmatrix}0&1\\1&0\end{smallmatrix}\bigr)$. If $R$ is
invertible, then it maps every summand of~\eqref{eq:charge-decomposition} injectively, hence
bijectively, onto itself. Consequently every $R_{ij}$ is invertible and $R^{-1}$ is charge-conserving.
Moreover $R(e_i\otimes e_j)=R_{ij}^{ij}\,e_i\otimes e_j+R_{ij}^{ji}\,e_j\otimes e_i$ is nonzero, and therefore
\begin{equation}\label{eq:line-criterion}
 R(L_i\otimes L_j)=L_j\otimes L_i\iff R_{ij}^{ij}=0 .
\end{equation}
For $\varnothing\neq J\subseteq I$, a charge-conserving $R$ preserves $V_J^{\otimes2}$, and
if $R\in\YBCC(V,\cL)$ then $R|_{V_J^{\otimes2}}\in\YBCC(V_J,\cL_J)$.

\begin{example}\label{ex:hecke}
Let $I=\{1,\dots,N\}$ and $Q\in\C^\times$ with $Q^2\neq1$. The \emph{standard Hecke operator}
is
\[
 R(e_i\otimes e_i)=Q\,e_i\otimes e_i,\quad
 R(e_i\otimes e_j)=(Q-Q^{-1})\,e_i\otimes e_j+e_j\otimes e_i,\quad
 R(e_j\otimes e_i)=e_i\otimes e_j\qquad(i<j).
\]
It is the operator of \cite[Section~4.2]{Tur} for the value $-Q$ of Turaev's parameter $q$
(see also \cite{Jim}); it is
charge-conserving, with $v_i=Q$ and $R_{ij}=\bigl(\begin{smallmatrix}Q-Q^{-1}&1\\1&0
\end{smallmatrix}\bigr)$ for $i<j$. The \emph{scaled flip} $\lambda\tau$,
$\tau(x\otimes y)=y\otimes x$, $\lambda\in\C^\times$, is charge-conserving for every
indexed line decomposition, with $v_i=\lambda$ and $R_{ij}=\lambda\tau|_{W_{ij}}$. In the
terminology of Section~\ref{subsec:components}, the Hecke operator with $N\ge2$ is the normal
form with scalar classes $\{1\}\prec\dots\prec\{N\}$, class scalars all equal to $Q$ and
spectral parameters $Q,-Q^{-1}$, and the scaled flip has the singletons as connected
components, each with $a=b=\lambda$, and $\kappa_{st}=\lambda^2$.
\end{example}

\subsection{Relabellings and X-transformations}
\label{subsec:PX}
A permutation $\varphi$ of $I$ acts by $\varphi\cdot F=(P_\varphi\otimes P_\varphi)F(P_\varphi\otimes
P_\varphi)^{-1}$, where $P_\varphi e_i:=e_{\varphi(i)}$; thus
$(\varphi\cdot F)_{\varphi(i)\varphi(j)}^{\varphi(k)\varphi(l)}=F_{ij}^{kl}$, and
$\rho_{\varphi\cdot R,r}(\xi)=P_\varphi^{\otimes r}\rho_{R,r}(\xi)(P_\varphi^{\otimes r})^{-1}$.
Let $\Diag(V^{\otimes2})$ be the group of operators acting by a scalar $y_{ij}\in\C^\times$ on
every line $L_i\otimes L_j$. For $Y\in\Diag(V^{\otimes2})$ one has
$(YFY^{-1})_{ij}^{kl}=(y_{kl}/y_{ij})F_{ij}^{kl}$. Hence $YFY^{-1}$ is charge-conserving when
$F$ is, its coefficients $F_{ij}^{ij}$ are unchanged, and
\begin{equation}\label{eq:x-action}
 (YFY^{-1})_{ij}=\begin{pmatrix}F_{ij}^{ij}&x_{ij}F_{ji}^{ij}\\[2pt]
 x_{ij}^{-1}F_{ij}^{ji}&F_{ji}^{ji}\end{pmatrix},\qquad x_{ij}:=\frac{y_{ij}}{y_{ji}} .
\end{equation}
Conjugation by such a $Y$ is called an \emph{X-transformation} \cite[Lemma~3.15]{MR}. It is a
change of basis of $V\otimes V$ that is diagonal in the basis $(e_i\otimes e_j)$, but in general
it is not induced by a change of basis of $V$, and therefore Lemma~\ref{lem:conjugation} does not apply to
it. It rescales the off-diagonal entries of each block and leaves its diagonal entries
unchanged. Every family $(x_{ij})_{i\ne j}$
with $x_{ji}=x_{ij}^{-1}$ arises in this way (take $y_{ij}:=x_{ij}$ for $i<j$ in some order of $I$, and $y_{ij}:=1$ otherwise).

Two charge-conserving operators $R,R'$ are \emph{PX-equivalent} if $R'=Y(\varphi\cdot R)Y^{-1}$
for some permutation $\varphi$ of $I$ and some $Y\in\Diag(V^{\otimes2})$, and
\emph{X-equivalent} if this holds with $\varphi=\id$. By the formula for $\rho_{\varphi\cdot R,r}$
and \cite[Lemma~3.15]{MR}, PX-equivalent operators are simultaneously Yang--Baxter operators,
and by a result of Martin and Rowell~\cite{MR2}, recalled in Proposition~\ref{prop:x-conjugacy}
below, their braid group representations are conjugate.

\subsection{The structure theorem of Martin and Rowell}
\label{subsec:components}
Following~\cite{MR}, we think of $R\in\YBCC(V,\cL)$ as the complete graph on $I$ whose vertex
$i$ carries the vertex scalar $v_i$ and whose edge $\{i,j\}$ carries the block $R_{ij}$. The
edge $\{i,j\}$ is \emph{scalar} if $R_{ij}=\lambda\,\id$ for some $\lambda$, a \emph{flip edge}
if $R_{ij}^{ij}=R_{ji}^{ji}=0$, and \emph{triangular} otherwise. The classification of the
two-dimensional case in~\cite{MR}, applied to the restriction of $R$ to
$V_{\{i,j\}}^{\otimes2}$, describes a single edge.

\begin{proposition}[{\cite[Proposition~3.22]{MR}}]\label{prop:two-index}
Let $R\in\YBCC(V,\cL)$ and $i\neq j$. Up to an X-transformation, the block $R_{ij}$ is one of
the following four matrices, where $\lambda,\eta,a,b\in\C^\times$ and $a+b\neq0$.
\begin{enumerate}[label=\textup{(\roman*)}]
\item $R_{ij}=\lambda\,\id$, and then $v_i=v_j=\lambda$.
\item $R_{ij}=\begin{pmatrix}0&\eta\\ \eta&0\end{pmatrix}$, that is, $R_{ij}$ is $\eta$ times the
flip of $W_{ij}$.
\item $R_{ij}=\begin{pmatrix}a+b&-ab\\1&0\end{pmatrix}$, and then $v_i,v_j\in\{a,b\}$.
\item $R_{ij}=\begin{pmatrix}0&-ab\\1&a+b\end{pmatrix}$, and then $v_i,v_j\in\{a,b\}$.
\end{enumerate}
\end{proposition}

Since X-transformations preserve the diagonal entries and the zero pattern of $R_{ij}$, the
edge $\{i,j\}$ is scalar in case~(i), a flip edge in case~(ii), and triangular in cases~(iii)
and~(iv). In case~(iii) we have $R(L_j\otimes L_i)=L_i\otimes L_j$ and
$R(L_i\otimes L_j)\neq L_j\otimes L_i$, in case~(iv) the reverse holds, and in both cases the
eigenvalues of $R_{ij}$ are $a$ and $b$. In~\cite{MR} the four cases are called $0$, $/$, $+$
and $-$, respectively. The Yang--Baxter equation on triples of indices
then forces the following global structure, which is the core of the classification
in~\cite{MR}.

\begin{proposition}[Martin--Rowell]\label{prop:MR-structure}
Let $R\in\YBCC(V,\cL)$.
\begin{enumerate}[label=\textup{(\roman*)}]
\item The relation $\sim_R$ on $I$, defined by $i\sim_Rj$ if $i=j$ or $\{i,j\}$ is not a flip edge, is an equivalence relation \cite[Lemma~6.6(I)]{MR}.
\item If $\{i,j\}$ and $\{i',j'\}$ are flip edges with $i\sim_Ri'$ and $j\sim_Rj'$, then
$R_{ji}^{ij}R_{ij}^{ji}=R_{j'i'}^{i'j'}R_{i'j'}^{j'i'}$ \cite[Lemma~6.6(II)]{MR}.
\item Suppose that $\sim_R$ has only one class. Then the relation $\approx_R$ on $I$, defined by $i\approx_Rj$ if $i=j$ or $\{i,j\}$ is a scalar edge, is an equivalence relation. If $C\neq D$ are classes of
$\approx_R$, then either $R(L_j\otimes L_i)=L_i\otimes L_j$ for all $i\in C$ and $j\in D$, or
$R(L_i\otimes L_j)=L_j\otimes L_i$ for all $i\in C$ and $j\in D$; writing $C\prec D$ in the first case
defines a total order on the set of classes \cite[Lemma~6.7(I)]{MR}.
\item In the situation of \textup{(iii)}, if there are at least two classes, then there are
$a,b\in\C^\times$ with $a+b\neq0$, unique up to order, such that every triangular block $R_{ij}$
has eigenvalues $a,b$, and $v_i\in\{a,b\}$ for every $i\in I$
\cite[Proposition~3.22 and Lemma~6.8(II)]{MR}.
\end{enumerate}
\end{proposition}

We call the classes $\Lambda_1,\dots,\Lambda_m$ of $\sim_R$ the \emph{connected components} of $R$ (the nations of~\cite{MR}), and we say that $R$ is \emph{connected} if $m=1$. We call $R^{(s)}:=R|_{V_{\Lambda_s}^{\otimes2}}$ the \emph{connected restriction} of $R$ to $\Lambda_s$; it is connected. For $i\in\Lambda_s$ and $j\in\Lambda_t$
with $s\neq t$, the edge $\{i,j\}$ is a flip edge, hence $R_{ij}^2=R_{ji}^{ij}R_{ij}^{ji}\,\id$; by
Proposition~\ref{prop:MR-structure}(ii) the scalar \[ \kappa_{st}:=R_{ji}^{ij}R_{ij}^{ji}\in\C^\times \] depends only on $s$ and $t$, and $\kappa_{st}=\kappa_{ts}$. By~\eqref{eq:x-action}, the connected components and the scalars
$\kappa_{st}$ are invariant under X-transformations.

Now let $R$ be connected. We call the classes of the relation $\approx_R$ of
Proposition~\ref{prop:MR-structure}(iii) the \emph{scalar classes} of $R$ (the counties of~\cite{MR}). If $C$ is a scalar class, then $R$ acts on $V_C^{\otimes2}$ by a scalar $v(C)$,
because for $i\neq j$ in $C$ the block $R_{ij}=\lambda\,\id$ is fixed by all X-transformations, and therefore Proposition~\ref{prop:two-index} gives $v_i=v_j=\lambda$. We label
the scalar classes $C_1\prec\dots\prec C_\ell$ and call $v_j:=v(C_j)$ the \emph{class scalars} of $R$. The vertex scalar of every $i\in C_j$ is therefore the class scalar $v_j$. Below, an index $j$ with $1\le j\le\ell$ on the letter $v$ refers to a class, and we say explicitly when a vertex scalar is meant.
If $\ell\ge2$, we call the scalars $a,b$ of Proposition~\ref{prop:MR-structure}(iv) the \emph{spectral parameters} of $R$; if $\ell=1$ the spectral parameters are $a:=v_1$ and $b:=v_1$. In both cases
$a+b\neq0$ and $v_j\in\{a,b\}$ for all $j$.

Martin and Rowell also construct all solutions explicitly. We need their construction for
connected operators and for the following way of combining them.

\begin{definition}\label{def:braided-sum}
Let $I=\Gamma_1\sqcup\dots\sqcup\Gamma_m$ with $\Gamma_s\neq\varnothing$, let
$S^{(s)}\in\EndCC(V_{\Gamma_s}^{\otimes2};\cL_{\Gamma_s})$ be invertible, and for $s\neq t$ let $\eta_{st}\in\C^\times$ with $\eta_{st}=\eta_{ts}$. The \emph{braided direct sum}
$S^{(1)}\bsum_\eta\dotsb\bsum_\eta S^{(m)}$ is the operator on $V^{\otimes2}$ that acts by $S^{(s)}$ on
$V_{\Gamma_s}^{\otimes2}$ and by $\eta_{st}\tau$ on
$(V_{\Gamma_s}\otimes V_{\Gamma_t})\oplus(V_{\Gamma_t}\otimes V_{\Gamma_s})$ for $s\neq t$.
\end{definition}

\begin{theorem}[Martin--Rowell]\label{thm:normal-form}
Let $I=C_1\sqcup\dots\sqcup C_\ell$ be an ordered partition into nonempty sets, let
$a,b\in\C^\times$ with $a+b\neq0$, with $a=b$ if $\ell=1$, and let
$v_1,\dots,v_\ell\in\{a,b\}$. Let $R^{\mathrm{nf}}\in\EndCC(V^{\otimes2};\cL)$ be the operator
given by
\begin{equation}\label{eq:normal-form}
 R^{\mathrm{nf}}|_{V_{C_r}^{\otimes2}}=v_r\,\id,\qquad
 R^{\mathrm{nf}}_{ij}=\begin{pmatrix}a+b&-ab\\1&0\end{pmatrix}\quad(i\in C_r,\ j\in C_{r'},\
 r<r'),
\end{equation}
the block being written in the basis $(e_i\otimes e_j,e_j\otimes e_i)$.
\begin{enumerate}[label=\textup{(\alph*)}]
\item $R^{\mathrm{nf}}$ is a connected Yang--Baxter operator with scalar classes
$C_1\prec\dots\prec C_\ell$, class scalars $v_1,\dots,v_\ell$ and, if $\ell\ge2$, spectral
parameters $a,b$.
\item Let $I=\Gamma_1\sqcup\dots\sqcup\Gamma_m$ with $\Gamma_s\neq\varnothing$, and for every
$s$ let $S^{(s)}$ be an operator of the form~\eqref{eq:normal-form} on
$(V_{\Gamma_s},\cL_{\Gamma_s})$. Then, for all $\eta_{st}\in\C^\times$ with $\eta_{st}=\eta_{ts}$, the braided
direct sum $R=S^{(1)}\bsum_\eta\dotsb\bsum_\eta S^{(m)}$ is a Yang--Baxter operator, its
connected components are $\Gamma_1,\dots,\Gamma_m$, $R^{(s)}=S^{(s)}$ and
$\kappa_{st}=\eta_{st}^2$.
\end{enumerate}
\end{theorem}

We call~\eqref{eq:normal-form} the \emph{normal form} with the given data. The Yang--Baxter
equation in (a) and (b) is \cite[Lemmas~6.9--6.11]{MR}. Indeed, up to X-transformations, which preserve the Yang--Baxter
equation \cite[Lemma~3.15]{MR}, the normal forms and their braided direct sums are the operators
produced by the recipe \cite[(4.10)]{MR}; the blocks of types $+$ and $-$ in
\cite[Proposition~3.22]{MR} are, up to X-transformations, the transposes of the block
in~\eqref{eq:normal-form} for $i<j$ and for $i>j$ in the order of $I=\{1,\dots,N\}$ used in~\cite{MR}, respectively. The other assertions of (a)
and (b) follow from the definitions.

\subsection{Two consequences}
\label{subsec:connected}
First, the off-diagonal entries of the flip edges between two components only matter through
their product $\kappa_{st}$, and therefore they can be normalized simultaneously.

\begin{proposition}\label{thm:canonical-decomposition}
Let $R\in\YBCC(V,\cL)$ have connected components $\Lambda_1,\dots,\Lambda_m$, and choose $\eta_{st}\in\C^\times$ with $\eta_{st}=\eta_{ts}$ and $\eta_{st}^2=\kappa_{st}$. Then $R$ is X-equivalent to
$R^{(1)}\bsum_\eta\dotsb\bsum_\eta R^{(m)}$ by an X-transformation that fixes $R$ on every
$V_{\Lambda_s}^{\otimes2}$, and different choices of the $\eta_{st}$ give X-equivalent
operators.
\end{proposition}

\begin{proof}
This is the normalization of the flip edges in the proof of \cite[Lemma~6.6(II)]{MR}; we only
make the X-transformation explicit. For $i\in\Lambda_s$ and $j\in\Lambda_t$ with $s\neq t$ put $x_{ij}:=\eta_{st}/R_{ji}^{ij}$, and put $x_{ij}:=1$ if $i$ and $j$ lie in the same component.
Then $x_{ji}=x_{ij}^{-1}$, because $R_{ji}^{ij}R_{ij}^{ji}=\eta_{st}^2$, and
by~\eqref{eq:x-action} the corresponding X-transformation fixes $R$ on every
$V_{\Lambda_s}^{\otimes2}$ and replaces every block between different components by
$\bigl(\begin{smallmatrix}0&\eta_{st}\\ \eta_{st}&0\end{smallmatrix}\bigr)$, that is, by
$\eta_{st}\tau$. If $\eta'_{st}=\pm\eta_{st}$ is another choice, the X-transformation with $x_{ij}:=\eta'_{st}/\eta_{st}$ for $i\in\Lambda_s$, $j\in\Lambda_t$, $s\neq t$, and $x_{ij}:=1$ otherwise, maps one braided direct sum to the other.
\end{proof}

Second, a connected operator satisfies a quadratic relation, because each of its blocks does.

\begin{corollary}\label{cor:quadratic}
If $R$ is connected with spectral parameters $a,b$, then $(R-a\,\id)(R-b\,\id)=0$ and
$R^{-1}=(ab)^{-1}\bigl((a+b)\id-R\bigr)$.
\end{corollary}

\begin{proof}
By~\eqref{eq:charge-decomposition}, $V^{\otimes2}$ is the sum of the lines $L_i\otimes L_i$ and
of the planes $W_{ij}$, and $R$ preserves all of them. On $L_i\otimes L_i$ it acts by the vertex scalar $v_i\in\{a,b\}$. On $W_{ij}$ it acts by the block $R_{ij}$, which is $v_i\,\id$ if $\{i,j\}$ is a scalar edge (Proposition~\ref{prop:two-index}), and has the
eigenvalues $a$ and $b$ otherwise (Proposition~\ref{prop:MR-structure}(iv)), since there are no
flip edges. In both cases $(R_{ij}-a)(R_{ij}-b)=0$, by the Cayley--Hamilton theorem in the
second case. This proves the first identity, and the second follows from
$R\bigl((a+b)\id-R\bigr)=ab\,\id$ (compare \cite[(A.2)]{MR}).
\end{proof}
\section{Diagonal enhancements}
\label{sec:enhancements}

In this section we determine all diagonal enhancements. After the definition, we show that
X-transformations conjugate the braid group representations and therefore do not change the invariants (Proposition~\ref{prop:px-invariance}). We then rewrite the enhancement equations in
terms of the diagonal coefficients $R_{kl}^{kl}$ and $(R^{-1})_{kl}^{kl}$
(Lemma~\ref{lem:coefficient-form}), and show that they split into one system for each connected
component (Proposition~\ref{prop:decoupling}). For a connected operator the system is triangular
with respect to the order of the scalar classes (Proposition~\ref{prop:connected-equations}),
and we solve it explicitly (Theorem~\ref{thm:connected-existence}).

\subsection{Definition and first properties}

\begin{definition}\label{def:enhancement}
Let $R\in\YBCC(V,\cL)$. An operator $\mu\in\End(V)$ is \emph{diagonal} if it preserves every line $L_i$. A \emph{diagonal enhancement} of $R$ is a triple $(R,\mu,\alpha)$, where $\mu\in\GL(V)$ is diagonal and $\alpha\in\C^\times$, such that
\begin{equation}\label{eq:enhancement}
 \tr_2\bigl(R(\id\otimes\mu)\bigr)=\alpha\,\id,\qquad
 \tr_2\bigl(R^{-1}(\id\otimes\mu)\bigr)=\alpha^{-1}\id .
\end{equation}
We write $\mu_i$ for the eigenvalue of $\mu$ on $L_i$, and set
\begin{equation}\label{eq:invariant}
 T_{R,\mu,\alpha}(\widehat\xi):=\alpha^{-w(\xi)}\tr\bigl(\rho_{R,r}(\xi)\mu^{\otimes r}\bigr) \qquad(\xi\in B_r).
\end{equation}
\end{definition}

\begin{lemma}\label{lem:turaev-relation}
Let $\mu\in\GL(V)$ preserve every $L_i$, and let $\alpha,\beta\in\C^\times$. Then
$(R,\mu,\alpha)$ is a diagonal enhancement if and only if $(R,\beta\mu,\alpha,\beta)$ is an
enhanced Yang--Baxter operator, and then $T_{(R,\beta\mu,\alpha,\beta)}=T_{R,\mu,\alpha}$. In
particular, $T_{R,\mu,\alpha}$ is a well-defined invariant of oriented links,
$T_{R,\mu,\alpha}(\bigcirc_r)=\tr(\mu)^r$, and $T_{R,\mu,\alpha}$ is multiplicative under
split union.
\end{lemma}

\begin{proof}
Condition (i) holds because $\mu\otimes\mu$ acts by a scalar on every summand
of~\eqref{eq:charge-decomposition}, namely by $\mu_i^2$ on $L_i\otimes L_i$ and by $\mu_i\mu_j$ on
$W_{ij}$, and $R$ preserves these summands. Since $\beta\mu$ is invertible,
condition (ii) for $(R,\beta\mu,\alpha,\beta)$ is equivalent to
$\tr_2\bigl(R^{\pm1}(\id\otimes\beta\mu)\bigr)=\alpha^{\pm1}\beta\,\id$, which
is~\eqref{eq:enhancement} multiplied by $\beta$. The factor $\beta^r$ of
$(\beta\mu)^{\otimes r}$ cancels $\beta^{-r}$, hence the two invariants agree, and the last
assertions follow from Theorem~\ref{thm:turaev}.
\end{proof}

\begin{lemma}\label{lem:coefficient-form}
Let $F\in\EndCC(V^{\otimes2};\cL)$ and let $\mu$ be diagonal. Then
$\tr_2\bigl(F(\id\otimes\mu)\bigr)$ acts on $L_k$ by $\sum_{l\in I}F_{kl}^{kl}\mu_l$.
Consequently $(R,\mu,\alpha)$ is a diagonal enhancement if and only if, for every $k\in I$,
\begin{equation}\label{eq:coefficient-form}
 \sum_{l\in I}R_{kl}^{kl}\,\mu_l=\alpha,\qquad\sum_{l\in I}(R^{-1})_{kl}^{kl}\,\mu_l=\alpha^{-1}.
\end{equation}
\end{lemma}

\begin{proof}
We compute the partial trace in the basis $(e_l)$, with dual basis $(e^l)$. Since
$\mu e_l=\mu_le_l$ and $F(e_k\otimes e_l)=\sum_{k',l'}F_{kl}^{k'l'}\,e_{k'}\otimes e_{l'}$,
\[
 \tr_2\bigl(F(\id\otimes\mu)\bigr)(e_k)=\sum_l\mu_l\,(\id\otimes e^l)F(e_k\otimes e_l)
 =\sum_{l,k'}\mu_l\,F_{kl}^{k'l}\,e_{k'} .
\]
Because $F$ is charge-conserving, $F_{kl}^{k'l}=0$ unless $(k',l)\in\{(k,l),(l,k)\}$, and in
both cases $k'=k$. Hence $\tr_2\bigl(F(\id\otimes\mu)\bigr)(e_k)=\bigl(\sum_lF_{kl}^{kl}\mu_l\bigr)e_k$.
Applying this to $F=R$ and to $F=R^{-1}$, which is also charge-conserving, shows
that~\eqref{eq:enhancement} is equivalent to~\eqref{eq:coefficient-form}.
\end{proof}

The next proposition, due to Martin and Rowell, shows that X-transformations do not affect the
braid group representations, up to conjugation by a diagonal operator.

\begin{proposition}[{\cite[Theorem~3.4]{MR2}}]\label{prop:x-conjugacy}
Let $R\in\YBCC(V,\cL)$, $Y\in\Diag(V^{\otimes2})$ and $R':=YRY^{-1}$. For $r\ge1$ let
$D_r\in\GL(V^{\otimes r})$ act on $L_{i_1}\otimes\dots\otimes L_{i_r}$ by
$\prod_{p<p'}y_{i_pi_{p'}}$. Then $R'_k=D_rR_kD_r^{-1}$ for $1\le k<r$, and therefore
$\rho_{R',r}(\xi)=D_r\rho_{R,r}(\xi)D_r^{-1}$ for all $\xi\in B_r$.
\end{proposition}

In~\cite{MR2} the result is stated as the $\infty$-equivalence of $R$ and $R'$, that is, the
equivalence of their braid group representations on $V^{\otimes r}$ for every $r$. In a basis
$(e_i)$ with $e_i\in L_i$, the intertwiner constructed in its proof is the operator $D_r$ above,
and Proposition~\ref{prop:px-invariance} uses that it is diagonal.

\begin{proposition}\label{prop:px-invariance}
Let $R\in\YBCC(V,\cL)$, let $\mu\in\GL(V)$ be diagonal and $\alpha\in\C^\times$. Let $R':=Y(\varphi\cdot R)Y^{-1}$ with $Y\in\Diag(V^{\otimes2})$ and $\varphi$ a permutation of $I$, and let $\mu':=P_\varphi\mu P_\varphi^{-1}$, that is, $\mu'_{\varphi(i)}=\mu_i$. Then
$(R,\mu,\alpha)$ is a diagonal enhancement if and only if $(R',\mu',\alpha)$ is one, and in
that case $T_{R',\mu',\alpha}=T_{R,\mu,\alpha}$. Moreover, the coefficients $R_{kl}^{kl}$ and $(R^{-1})_{kl}^{kl}$ are invariant under X-transformations.
\end{proposition}

\begin{proof}
Conjugation by $Y$ does not change the coefficients $F_{kl}^{kl}$, and
$(YFY^{-1})^{-1}=YF^{-1}Y^{-1}$; moreover $(\varphi\cdot R)^{-1}=\varphi\cdot R^{-1}$ and
$(\varphi\cdot F)_{\varphi(k)\varphi(l)}^{\varphi(k)\varphi(l)}=F_{kl}^{kl}$. Hence
\[
 \sum_{l\in I}(R'^{\pm1})_{\varphi(k)\varphi(l)}^{\varphi(k)\varphi(l)}\,\mu'_{\varphi(l)}
 =\sum_{l\in I}(R^{\pm1})_{kl}^{kl}\,\mu_l,
\] and therefore the equations~\eqref{eq:coefficient-form} for $(R',\mu',\alpha)$ at $\varphi(k)$ coincide with
those for $(R,\mu,\alpha)$ at $k$. By Proposition~\ref{prop:x-conjugacy},
$\rho_{R',r}(\xi)=D_rP_\varphi^{\otimes r}\rho_{R,r}(\xi)(P_\varphi^{\otimes r})^{-1}D_r^{-1}$, where
$D_r$ commutes with $\mu'^{\otimes r}$ because both act by scalars on every
$L_{i_1}\otimes\dots\otimes L_{i_r}$, and
$\mu'^{\otimes r}=P_\varphi^{\otimes r}\mu^{\otimes r}(P_\varphi^{\otimes r})^{-1}$. The equality of
the invariants follows from the cyclicity of the trace.
\end{proof}

\subsection{Reduction to connected components}
For the rest of the paper, $R\in\YBCC(V,\cL)$ has connected components
$\Lambda_1,\dots,\Lambda_m$ and connected restrictions $R^{(s)}$; the scalar classes of $R^{(s)}$
are $C_{s,1}\prec\dots\prec C_{s,\ell_s}$, with class scalars $v_{s,j}$ and a chosen ordering
$(a_s,b_s)$ of the spectral parameters.

\begin{lemma}\label{lem:diagonal-coefficients}
Let $k\in C_{s,j}$ and $l\in I$. The coefficients $R_{kl}^{kl}$ and $(R^{-1})_{kl}^{kl}$ are
\begin{center}
\renewcommand{\arraystretch}{1.25}
\begin{tabular}{c|c|c}
position of $l$ & $R_{kl}^{kl}$ & $(R^{-1})_{kl}^{kl}$\\ \hline
$l\in C_{s,j}$ & $v_{s,j}$ & $v_{s,j}^{-1}$\\
$l\in C_{s,r}$, $r>j$ & $a_s+b_s$ & $0$\\
$l\in C_{s,r}$, $r<j$ & $0$ & $(a_s+b_s)/(a_sb_s)$\\
$l\notin\Lambda_s$ & $0$ & $0$
\end{tabular}
\end{center}
\end{lemma}

\begin{proof}
\emph{The last row.} If $l\notin\Lambda_s$, then $\{k,l\}$ is a flip edge. Hence $R$ maps $L_k\otimes L_l$ onto $L_l\otimes L_k$, and $R^{-1}$ does as well. Hence $R(e_k\otimes e_l)$ and
$R^{-1}(e_k\otimes e_l)$ have no component on $e_k\otimes e_l$, and both coefficients vanish.

\emph{The other rows.} If $l\in\Lambda_s$, the two coefficients are those of $R^{(s)}$ and of
$(R^{(s)})^{-1}$, because $R$ preserves $V_{\Lambda_s}^{\otimes2}$ and
$R^{-1}|_{V_{\Lambda_s}^{\otimes2}}=(R^{(s)})^{-1}$. If $l\in C_{s,j}$, then $R$ acts on
$L_k\otimes L_l$ by the class scalar $v_{s,j}$ (Section~\ref{subsec:components}). If
$l\in C_{s,r}$ with $r>j$, then $C_{s,j}\prec C_{s,r}$, hence $R(L_l\otimes L_k)=L_k\otimes L_l$
and $R_{lk}^{lk}=0$ by~\eqref{eq:line-criterion}. The block $R_{kl}$ has trace
$R_{kl}^{kl}+R_{lk}^{lk}$ and eigenvalues $a_s,b_s$ (Proposition~\ref{prop:MR-structure}(iv)),
hence $R_{kl}^{kl}=a_s+b_s$. If $l\in C_{s,r}$ with $r<j$, then $C_{s,r}\prec C_{s,j}$, hence $R(L_k\otimes L_l)=L_l\otimes L_k$ and $R_{kl}^{kl}=0$ by~\eqref{eq:line-criterion}. This
gives the column of $R_{kl}^{kl}$. For the
column of $(R^{-1})_{kl}^{kl}$ we use Corollary~\ref{cor:quadratic}, which gives
$R^{-1}=(a_sb_s)^{-1}\bigl((a_s+b_s)\id-R\bigr)$ on $V_{\Lambda_s}^{\otimes2}$, and therefore $(R^{-1})_{kl}^{kl}=(a_s+b_s-R_{kl}^{kl})/(a_sb_s)$. For $R_{kl}^{kl}=v\in\{a_s,b_s\}$ this is
$v^{-1}$, since $(a_s+b_s-v)v=a_sb_s$; for $R_{kl}^{kl}=a_s+b_s$ it is $0$; and for
$R_{kl}^{kl}=0$ it is $(a_s+b_s)/(a_sb_s)$.
\end{proof}

\begin{proposition}\label{prop:decoupling}
Let $\mu\in\GL(V)$ be diagonal and $\alpha\in\C^\times$, and put $\mu^{(s)}:=\mu|_{V_{\Lambda_s}}$. Then $(R,\mu,\alpha)$ is a diagonal enhancement if and only
if $(R^{(s)},\mu^{(s)},\alpha)$ is a diagonal enhancement for every $s$.
\end{proposition}

\begin{proof}
For $k\in\Lambda_s$, the terms with $l\notin\Lambda_s$ in~\eqref{eq:coefficient-form} vanish by
Lemma~\ref{lem:diagonal-coefficients}, and the remaining ones form the
equations~\eqref{eq:coefficient-form} for $(R^{(s)},\mu^{(s)},\alpha)$ at $k$.
\end{proof}

In particular, the scalars $\kappa_{st}$ play no role in the existence of diagonal
enhancements; they will enter the invariants.

\subsection{Connected operators}
In this subsection $R$ is connected, with scalar classes $C_1\prec\dots\prec C_\ell$, class
scalars $v_j$, and an ordering $(a,b)$ of its spectral parameters. Recall that $a+b\neq0$. Put $\bar v_j:=a+b-v_j$, which satisfies $v_j\bar v_j=ab$, and, for a diagonal $\mu$, let $t_j:=\sum_{k\in C_j}\mu_k$, the trace of $\mu$ on $V_{C_j}$.

\begin{proposition}\label{prop:connected-equations}
Let $\mu\in\GL(V)$ be diagonal and $\alpha\in\C^\times$. The following conditions are
equivalent.
\begin{enumerate}[label=\textup{(\roman*)}]
\item $(R,\mu,\alpha)$ is a diagonal enhancement;
\item $v_jt_j+(a+b)\sum_{r>j}t_r=\alpha$ for $1\le j\le\ell$, and
$(a+b)\tr(\mu)=\alpha+ab\,\alpha^{-1}$;
\item $t_1=v_1\alpha^{-1}$, $\ t_{j+1}=-\dfrac{v_jv_{j+1}}{ab}\,t_j$ for $1\le j<\ell$, and
$v_\ell t_\ell=\alpha$.
\end{enumerate}
\end{proposition}

\begin{proof}
\emph{Step 1. The first equation.} Put $P_j:=v_jt_j+(a+b)\sum_{r>j}t_r$, and let $k\in C_j$. By
Lemma~\ref{lem:coefficient-form}, $\tr_2\bigl(R(\id\otimes\mu)\bigr)$ acts on $L_k$ by
$\sum_lR_{kl}^{kl}\mu_l$. By Lemma~\ref{lem:diagonal-coefficients}, the indices $l$ in $C_j$
contribute $v_j\sum_{l\in C_j}\mu_l=v_jt_j$, the indices in a later class $C_r$ contribute
$(a+b)t_r$, and the indices in earlier classes contribute nothing. Hence
$\tr_2\bigl(R(\id\otimes\mu)\bigr)$ acts on $L_k$ by $P_j$.

\emph{Step 2. The second equation.} By Corollary~\ref{cor:quadratic},
$R^{-1}=(ab)^{-1}\bigl((a+b)\id-R\bigr)$, and $\tr_2(\id\otimes\mu)=\tr(\mu)\,\id$. Therefore
\[
 \tr_2\bigl(R^{-1}(\id\otimes\mu)\bigr)
 =(ab)^{-1}\Bigl((a+b)\tr(\mu)\,\id-\tr_2\bigl(R(\id\otimes\mu)\bigr)\Bigr).
\]
Consequently~\eqref{eq:enhancement} holds if and only if $P_j=\alpha$ for all $j$ and
$(ab)^{-1}\bigl((a+b)\tr(\mu)-\alpha\bigr)=\alpha^{-1}$, that is,
$(a+b)\tr(\mu)=\alpha+ab\,\alpha^{-1}$. This proves (i)$\Leftrightarrow$(ii).

\emph{Step 3. The recursion.} The conditions $P_1=\dots=P_\ell=\alpha$ form a triangular system,
since $P_\ell$ involves only $t_\ell$ and $P_j$ differs from $P_{j+1}$ only through $t_j$ and
$t_{j+1}$. Indeed,
\[
 P_j-P_{j+1}=v_jt_j+\bar v_{j+1}t_{j+1}\quad(1\le j<\ell),\qquad P_\ell=v_\ell t_\ell,\qquad
 P_1=(a+b)\tr(\mu)-\bar v_1t_1 .
\]
The conditions $P_1=\dots=P_\ell=\alpha$ are equivalent to $P_\ell=\alpha$ together with
$P_j-P_{j+1}=0$ for $1\le j<\ell$, that is, to $v_\ell t_\ell=\alpha$ together with
$v_jt_j+\bar v_{j+1}t_{j+1}=0$; since $\bar v_{j+1}=ab/v_{j+1}$, the latter is the recursion
$t_{j+1}=-v_jv_{j+1}t_j/(ab)$ of (iii). Finally, when these conditions hold we have $P_1=\alpha$, hence the condition $(a+b)\tr(\mu)=\alpha+ab\,\alpha^{-1}$ of (ii) becomes
$\alpha+\bar v_1t_1=\alpha+ab\,\alpha^{-1}$, that is, $\bar v_1t_1=ab\,\alpha^{-1}$, or
$t_1=v_1\alpha^{-1}$. This proves (ii)$\Leftrightarrow$(iii).
\end{proof}

\subsection{Existence and parametrization}

\begin{definition}\label{def:signed-rank}
Let $R$ be connected as above. Put $\theta_j:=1$ if $v_j=a$, and $\theta_j:=-1$ if $v_j=b\neq a$; let $h_0:=0$, $h_j:=\theta_1+\dots+\theta_j$, and $q:=-a/b$. We define the \emph{signed effective rank} of $R$ to be $n:=h_\ell$, and we set \[ \cE(R):=-ab\,q^n=-ab\Bigl(-\frac ab\Bigr)^n .
\]
\end{definition}

If $a\neq b$, then $n$ is the number of classes with class scalar $a$ minus the number of
classes with class scalar $b$; if $a=b$, then $\theta_j=1$ for all $j$, $n=\ell$ and $q=-1$.

\begin{lemma}\label{lem:v-identities}
For $1\le j\le\ell$,
\begin{equation}\label{eq:v-identities}
 -\frac{v_j^2}{ab}=q^{\theta_j},\qquad v_j=\frac{ab}{a+b}\bigl(1-q^{\theta_j}\bigr).
\end{equation}
Consequently
\begin{equation}\label{eq:E-product}
 \cE(R)=\frac{(-1)^{\ell-1}\bigl(\prod_{j=1}^{\ell}v_j\bigr)^2}{(ab)^{\ell-1}}\neq0,
\end{equation}
and $\cE(R)$ does not depend on the ordering of $a,b$.
\end{lemma}

\begin{proof}
If $v_j=a$ (including the case $a=b$), then $\theta_j=1$, and $-a^2/(ab)=-a/b=q$ and
$\frac{ab}{a+b}(1+\frac ab)=a$. If $v_j=b\neq a$, then $\theta_j=-1$, and
$-b^2/(ab)=-b/a=q^{-1}$ and $\frac{ab}{a+b}(1+\frac ba)=b$. Multiplying the first identity over
$j$ gives $\prod_jv_j^2=(-ab)^\ell q^{n}$, which is~\eqref{eq:E-product}. The right-hand side
of~\eqref{eq:E-product} is symmetric in $a,b$.
\end{proof}

\begin{theorem}\label{thm:connected-existence}
Let $R$ be connected and $\alpha\in\C^\times$. Then $R$ admits a diagonal enhancement
$(R,\mu,\alpha)$ if and only if $\alpha^2=\cE(R)$. In that case, a diagonal $\mu\in\GL(V)$
gives a diagonal enhancement $(R,\mu,\alpha)$ if and only if
\begin{equation}\label{eq:class-traces}
 t_j=\alpha^{-1}v_jq^{h_{j-1}}
 =\frac{ab}{\alpha(a+b)}\bigl(q^{h_{j-1}}-q^{h_j}\bigr)\qquad(1\le j\le\ell),
\end{equation}
and then $\tr(\mu)=\dfrac{ab}{\alpha(a+b)}\bigl(1-q^n\bigr)=\dfrac{\alpha^2+ab}{\alpha(a+b)}$.
\end{theorem}

\begin{proof}
By Proposition~\ref{prop:connected-equations}, we must solve the conditions (iii) there, namely
the initial value $t_1=v_1\alpha^{-1}$, the recursion, and the final condition
$v_\ell t_\ell=\alpha$.

\emph{Step 1. Solving the recursion.} Dividing the recursion by $v_{j+1}$ and
using~\eqref{eq:v-identities} gives
\[
 \frac{t_{j+1}}{v_{j+1}}=-\frac{v_j^2}{ab}\cdot\frac{t_j}{v_j}=q^{\theta_j}\,\frac{t_j}{v_j}.
\]
Hence the recursion together with $t_1/v_1=\alpha^{-1}=\alpha^{-1}q^{h_0}$ is equivalent to
$t_j/v_j=\alpha^{-1}q^{h_{j-1}}$ for all $j$, which is the first equality
in~\eqref{eq:class-traces}. The second equality follows by substituting the second identity
in~\eqref{eq:v-identities} for $v_j$ and using $q^{h_{j-1}}q^{\theta_j}=q^{h_j}$.

\emph{Step 2. The final condition.} For the values~\eqref{eq:class-traces}, the first identity
in~\eqref{eq:v-identities} gives
\[
 v_\ell t_\ell=\alpha^{-1}v_\ell^2q^{h_{\ell-1}}=-ab\,\alpha^{-1}q^{\theta_\ell+h_{\ell-1}}
 =-ab\,\alpha^{-1}q^{n},
\] and therefore the condition $v_\ell t_\ell=\alpha$ is equivalent to $\alpha^2=-ab\,q^n=\cE(R)$. Together
with Step~1, this shows that $(R,\mu,\alpha)$ is a diagonal enhancement if and only if
$\alpha^2=\cE(R)$ and~\eqref{eq:class-traces} holds.

\emph{Step 3. Existence.} The conditions~\eqref{eq:class-traces} only prescribe the sum of the
eigenvalues of $\mu$ on each scalar class, and these sums are nonzero. Hence, if $\alpha^2=\cE(R)$, we may take, for instance, $\mu_k=t_j/|C_j|$ for $k\in C_j$, and a diagonal enhancement exists.

\emph{Step 4. The trace.} Summing the second expression in~\eqref{eq:class-traces} over $j$,
the sum telescopes to $\tr(\mu)=\frac{ab}{\alpha(a+b)}(q^{h_0}-q^{h_\ell})=
\frac{ab}{\alpha(a+b)}(1-q^n)$, and $ab(1-q^n)=ab+\alpha^2$ because $-abq^n=\alpha^2$.
\end{proof}

\begin{example}\label{ex:three-classes-enh}
Let $I=\{1,2,3\}$, let $a,b\in\C^\times$ with $a\neq b$ and $a+b\neq0$, and let $R$ be the
operator~\eqref{eq:normal-form} with scalar classes $C_1=\{1\}\prec C_2=\{2,3\}$ and class
scalars $v_1=a$, $v_2=b$, that is,
\[
 R(e_1\otimes e_1)=a\,e_1\otimes e_1,\qquad R(e_i\otimes e_j)=b\,e_i\otimes e_j\quad
 (i,j\in\{2,3\}),
\]
\[
 R(e_1\otimes e_j)=(a+b)\,e_1\otimes e_j+e_j\otimes e_1,\qquad
 R(e_j\otimes e_1)=-ab\,e_1\otimes e_j\qquad(j\in\{2,3\}).
\]
By Theorem~\ref{thm:normal-form}(a) it is a connected Yang--Baxter operator. We have $(\theta_1,\theta_2)=(1,-1)$, $h_0=0$, $h_1=1$, $h_2=0$, hence $n=0$ and
$\cE(R)=-ab$. Thus $\alpha=\pm\sqrt{-ab}$, and~\eqref{eq:class-traces} gives
\[
 t_1=\frac a\alpha,\qquad t_2=\frac{b}{\alpha}\Bigl(-\frac ab\Bigr)=-\frac a\alpha .
\]
Hence the diagonal enhancements with a given $\alpha$ are the $\mu=\operatorname{diag}(\mu_1,\mu_2,
\mu_3)$ with $\mu_1=a/\alpha$ and $\mu_2+\mu_3=-a/\alpha$, $\mu_2,\mu_3\neq0$; they form a line
with two points removed, and $\tr(\mu)=0$. The equations of
Proposition~\ref{prop:connected-equations}(ii) can be checked directly, since $v_1t_1+(a+b)t_2=-ab/\alpha=\alpha$ and
$v_2t_2=-ab/\alpha=\alpha$. By Corollary~\ref{cor:degenerate}(i) below, the invariants of these
enhancements vanish.
\end{example}

\begin{theorem}[Theorem~\ref{mainthm:A}]\label{thm:global-existence}
An operator $R\in\YBCC(V,\cL)$ admits a diagonal enhancement if and only if
\[
 \cE(R^{(1)})=\dots=\cE(R^{(m)}).
\]
In that case, every diagonal enhancement $(R,\mu,\alpha)$ satisfies $\alpha^2=\cE(R^{(1)})$. Fix
such an $\alpha$, and let $t_{s,j}(\alpha)$ be the value~\eqref{eq:class-traces} for the class
$C_{s,j}$ of $R^{(s)}$. Then the set
$\cM_{R,\alpha}$ of diagonal $\mu\in\GL(V)$ such that $(R,\mu,\alpha)$ is a diagonal
enhancement is
\[
 \cM_{R,\alpha}=\prod_{s=1}^m\prod_{j=1}^{\ell_s}\cH_{s,j}(\alpha),\qquad \cH_{s,j}(\alpha):=\Bigl\{(\mu_k)_{k\in C_{s,j}}\in(\C^\times)^{C_{s,j}}:
 \sum_{k\in C_{s,j}}\mu_k=t_{s,j}(\alpha)\Bigr\}.
\]
Each $\cH_{s,j}(\alpha)$ is nonempty, and it consists of the points of an affine hyperplane of
$\C^{C_{s,j}}$ whose coordinates are all nonzero. Thus $\cM_{R,\alpha}$ is obtained from an
affine subspace of $\C^I$, of dimension $\dim V$ minus the total number of scalar classes, by
removing the points with a vanishing coordinate, and $\mu$ is unique when all scalar classes are
singletons.
\end{theorem}

\begin{proof}
By Proposition~\ref{prop:decoupling}, $(R,\mu,\alpha)$ is a diagonal enhancement if and only if
$(R^{(s)},\mu^{(s)},\alpha)$ is a diagonal enhancement for every $s$, with the same $\alpha$
for all $s$. By Theorem~\ref{thm:connected-existence}, this requires $\alpha^2=\cE(R^{(s)})$ for
every $s$, which is possible if and only if all the $\cE(R^{(s)})$ coincide. If they do, and
$\alpha^2$ is their common value, Theorem~\ref{thm:connected-existence} says that
$(R^{(s)},\mu^{(s)},\alpha)$ is a diagonal enhancement if and only if the eigenvalues of $\mu$
on each class $C_{s,j}$ add up to $t_{s,j}(\alpha)$. Different classes involve disjoint sets of eigenvalues, and therefore the set of solutions is the product of the sets $\cH_{s,j}(\alpha)$. Each of
them is nonempty because $t_{s,j}(\alpha)\neq0$; it is the complement, in the affine hyperplane
$\sum_{k\in C_{s,j}}\mu_k=t_{s,j}(\alpha)$, of the coordinate hyperplanes $\mu_k=0$.
\end{proof}

For Turaev's quadruples, Lemma~\ref{lem:turaev-relation} shows that the enhanced Yang--Baxter
operators with diagonal invertible enhancing endomorphism are the $(R,\beta\mu,\alpha,\beta)$
with $\beta\in\C^\times$ arbitrary and $\mu\in\cM_{R,\alpha}$.

\begin{remark}\label{rem:singular-mu}
Theorem~\ref{thm:global-existence} also describes the enhanced Yang--Baxter operators
$S=(R,\mu,\alpha,\beta)$ with $R\in\YBCC(V,\cL)$ and $\mu$ diagonal but not invertible. Let
$\mu_i$ be the eigenvalue of $\mu$ on $L_i$ and put $J:=\{i\in I:\mu_i\neq0\}$.
Condition~(i) in the definition of an enhanced Yang--Baxter operator
(Section~\ref{subsec:turaev}) always holds, because $\mu\otimes\mu$ acts by a scalar on every
summand of~\eqref{eq:charge-decomposition}. Since
$\tr_2\bigl(F(\mu\otimes\mu)\bigr)=\tr_2\bigl(F(\id\otimes\mu)\bigr)\mu$,
Lemma~\ref{lem:coefficient-form} shows that condition~(ii) holds on $L_k$ for every $k\notin J$,
and that for $k\in J$ it reads $\sum_{l\in J}(R^{\pm1})_{kl}^{kl}\,\beta^{-1}\mu_l=\alpha^{\pm1}$.
If $J=\varnothing$, that is, $\mu=0$, then $S$ is an enhanced Yang--Baxter operator for all
$\alpha,\beta\in\C^\times$, and $T_S$ vanishes on every nonempty link. If $J\neq\varnothing$, let
$R^{(J)}:=R|_{V_J^{\otimes2}}$, which belongs to $\YBCC(V_J,\cL_J)$ by
Section~\ref{subsec:blocks}. The coefficients of $R^{(J)}$ and of $(R^{(J)})^{-1}$ are those of
$R$ and of $R^{-1}$ with indices in $J$. Hence $S$ is an enhanced Yang--Baxter operator if and
only if $(R^{(J)},\beta^{-1}\mu|_{V_J},\alpha)$ is a diagonal enhancement of $R^{(J)}$. In that
case $T_S=T_{R^{(J)},\beta^{-1}\mu|_{V_J},\alpha}$, because $\rho_{R,r}(\xi)$ preserves
$V_J^{\otimes r}$, on which it acts by $\rho_{R^{(J)},r}(\xi)$, and it also preserves the sum of
the lines $L_{i_1}\otimes\dots\otimes L_{i_r}$ with some $i_p\notin J$, on which
$\mu^{\otimes r}$ vanishes.
\end{remark}

\subsection{Examples}

\begin{example}\label{ex:hecke-enhancement}
For the standard Hecke operator (Example~\ref{ex:hecke}) we have $a=Q$, $b=-Q^{-1}$, $-ab=1$,
$-a/b=Q^2$ and $n=N$. Hence $\cE(R)=Q^{2N}$ and $\alpha=\pm Q^N$. All scalar classes are
singletons, and for $\alpha=Q^N$ formula~\eqref{eq:class-traces} gives
$\mu=\operatorname{diag}(Q^{1-N},Q^{3-N},\dots,Q^{N-1})$; for $\alpha=-Q^N$ one gets $-\mu$.
In \cite[Theorem~4.2.1]{Tur} this operator appears with Turaev's parameter $q$ equal to $-Q$
(Example~\ref{ex:hecke}), enhanced by $\mu_i=q^{2i-N-1}$, $\alpha=-q^N$ and $\beta=1$. Since
$(-Q)^{2i-N-1}=(-1)^{N+1}Q^{2i-N-1}$ and $-(-Q)^N=(-1)^{N+1}Q^N$, this is the enhancement
$(R,(-1)^{N+1}\mu,(-1)^{N+1}Q^N)$, one of the two found above.
\end{example}

\begin{example}\label{ex:flip-enhancement}
For a scaled flip $\lambda\tau$, every component is a singleton with $\ell=1$, $a=b=\lambda$
and $n=1$, hence $\cE(R^{(s)})=\lambda^2$ for every $s$. Therefore the diagonal enhancements are $(\lambda\tau,\id,\lambda)$
and $(\lambda\tau,-\id,-\lambda)$.
\end{example}

\section{Link invariants}
\label{sec:invariants}

We now compute the invariants of the diagonal enhancements found in
Section~\ref{sec:enhancements}. After two elementary normalizations
(Section~\ref{subsec:normalizations}), we treat connected operators, whose invariants are
specializations of the HOMFLY--PT polynomial (Section~\ref{subsec:connected-invariants}). We
then show that the invariant of an arbitrary operator is a state sum over colourings of the
components of the link (Section~\ref{subsec:state-sums}), and finally that every such state sum
occurs (Section~\ref{subsec:realization}).

\subsection{Normalizations}
\label{subsec:normalizations}

\begin{proposition}\label{prop:scaling}
Let $(R,\mu,\alpha)$ be a diagonal enhancement and $L$ a link.
\begin{enumerate}[label=\textup{(\roman*)}]
\item For $\lambda\in\C^\times$, $(\lambda R,\mu,\lambda\alpha)$ is a diagonal enhancement and
$T_{\lambda R,\mu,\lambda\alpha}(L)=T_{R,\mu,\alpha}(L)$.
\item $(R,-\mu,-\alpha)$ is a diagonal enhancement and
$T_{R,-\mu,-\alpha}(L)=(-1)^{|L|}T_{R,\mu,\alpha}(L)$.
\end{enumerate}
\end{proposition}

\begin{proof}
The equations~\eqref{eq:enhancement} for $(\lambda R,\mu,\lambda\alpha)$ are those for $(R,\mu,\alpha)$ multiplied by $\lambda^{\pm1}$, and the equations for $(R,-\mu,-\alpha)$ are those for $(R,\mu,\alpha)$ multiplied by $-1$. For $\xi\in B_r$ we have
$\rho_{\lambda R,r}(\xi)=\lambda^{w(\xi)}\rho_{R,r}(\xi)$, which gives (i), and the factor
in (ii) is $(-1)^{r+w(\xi)}=(-1)^{|\widehat\xi|}$ by~\eqref{eq:parity}.
\end{proof}

Thus every invariant $T_{R,\mu,\alpha}$ is also the invariant of the diagonal enhancement
$(\alpha^{-1}R,\mu,1)$. The operations of Proposition~\ref{prop:scaling} multiply every
$\kappa_{st}$ and $\alpha^2$ by the same factor, and therefore the ratios $\kappa_{st}/\alpha^2$ are
unchanged.

\subsection{Connected invariants}
\label{subsec:connected-invariants}
Throughout this subsection $R$ is connected. The invariant will be identified by the method of
Lemma~\ref{lem:skein-uniqueness}, by showing that it satisfies a HOMFLY--PT skein relation and
computing it on unlinks. It is convenient to rescale the spectral parameters first.

\begin{definition}\label{def:effective-data}
Let $(R,\mu,\alpha)$ be a diagonal enhancement of a connected $R$, with an ordering $(a,b)$
of its spectral parameters and signed effective rank $n$. Choose $d\in\C^\times$ with
$d^2=-ab$ and put $Q:=a/d$ and $\varepsilon:=\alpha/(dQ^n)$. We call $(Q,n,\varepsilon)$ the
\emph{effective data} of $(R,\mu,\alpha)$ (for these choices).
\end{definition}

\begin{lemma}\label{lem:effective-data}
In Definition~\ref{def:effective-data}, $a=dQ$, $b=-dQ^{-1}$, $Q^2=-a/b\neq1$,
$\varepsilon\in\{\pm1\}$, and $\tr(\mu)=\varepsilon[n]_Q$.
\end{lemma}

\begin{proof}
We have $b=-d^2/a=-dQ^{-1}$, hence $Q^2=-a/b$, and $Q^2=1$ would mean $a+b=0$. By
Theorem~\ref{thm:connected-existence}, $\alpha^2=-ab(-a/b)^n=d^2Q^{2n}$, hence $\varepsilon^2=1$. Using $ab=-d^2$, $a+b=d(Q-Q^{-1})$ and $\alpha=\varepsilon dQ^n$ in the
formula of Theorem~\ref{thm:connected-existence},
\[
 \tr(\mu)=\frac{ab\,(1-Q^{2n})}{\alpha(a+b)}
 =\frac{-d^2(1-Q^{2n})}{\varepsilon dQ^n\cdot d(Q-Q^{-1})}=\varepsilon[n]_Q. \qedhere
\]
\end{proof}

\begin{theorem}[Theorem~\ref{mainthm:B}]\label{thm:connected-invariant}
Let $(R,\mu,\alpha)$ be a diagonal enhancement of a connected $R$, with effective data
$(Q,n,\varepsilon)$. Then, for every oriented link $L$,
\[
 T_{R,\mu,\alpha}(L)=\varepsilon[n]_Q\,P_L\bigl(\varepsilon Q^n,\,Q-Q^{-1}\bigr).
\]
\end{theorem}

\begin{proof}
Let $T:=T_{R,\mu,\alpha}$ and put \[ A:=\frac\alpha d=\varepsilon Q^n,\qquad z:=\frac{a+b}d=Q-Q^{-1},\qquad \delta:=\tr(\mu)=\varepsilon[n]_Q,
\]
where the last equality is Lemma~\ref{lem:effective-data}. We compare $T$ with the invariant $H(L):=\delta P_L(A,z)$.

\emph{Step 1. The skein relation for $T$.} By Lemma~\ref{lem:turaev-relation}, $T$ is the
invariant $T_S$ of the enhanced Yang--Baxter operator $S=(R,\mu,\alpha,1)$, and $R$ satisfies
$(R-a)(R-b)=0$ by Corollary~\ref{cor:quadratic}. Hence Theorem~\ref{thm:turaev-skein} gives
\[
 \alpha\,T(L_+)-(a+b)\,T(L_0)+ab\,\alpha^{-1}\,T(L_-)=0
\]
for every skein triple. We divide by $d$. Since $\alpha/d=A$, $(a+b)/d=z$ and
$ab/(d\alpha)=-d^2/(d\alpha)=-A^{-1}$, this becomes $A\,T(L_+)-A^{-1}T(L_-)=z\,T(L_0)$, which
is the skein relation~\eqref{eq:homfly-intro} satisfied by $H$.

\emph{Step 2. Unlinks.} By Lemma~\ref{lem:turaev-relation}, $T(\bigcirc_r)=\tr(\mu)^r=\delta^r$.
On the other hand $H(\bigcirc_r)=\delta\bigl((A-A^{-1})/z\bigr)^{r-1}$, and
\[
 \frac{A-A^{-1}}{z}=\frac{\alpha/d-d/\alpha}{(a+b)/d}=\frac{\alpha^2-d^2}{\alpha(a+b)}
 =\frac{\alpha^2+ab}{\alpha(a+b)}=\delta
\]
by the last formula of Theorem~\ref{thm:connected-existence}. Hence
$H(\bigcirc_r)=\delta^r=T(\bigcirc_r)$.

\emph{Step 3. Conclusion.} By Lemma~\ref{lem:skein-uniqueness}, $T=H$. Note that the argument
also covers the case $\delta=0$, in which both invariants vanish identically.
\end{proof}

\begin{remark}\label{rem:choices}
The value in Theorem~\ref{thm:connected-invariant} does not depend on the choices. Replacing
$d$ by $-d$ changes $(Q,\varepsilon)$ into $(-Q,(-1)^{n+1}\varepsilon)$; this changes the
signs of $A$ and $z$ and fixes $\delta$, and therefore by~\eqref{eq:homfly-sign} the value is unchanged.
If $a\neq b$, exchanging $a$ and $b$ changes $(Q,n,\varepsilon)$ into
$(-Q^{-1},-n,(-1)^n\varepsilon)$, which fixes $A$, $z$ and $\delta$.
\end{remark}

\begin{remark}\label{rem:not-detected}
Two connected diagonal enhancements with the same effective data have the same invariant. In
particular, the invariant does not depend on $\dim V$, on the sizes of the scalar classes, on
the order of the class scalars beyond the integer $n$, or on the eigenvalues of $\mu$ inside
each scalar class. This is reminiscent of the fact that the $\mathfrak{gl}_{m|n}$ link invariants
depend only on $m-n$ \cite[Theorem~4.7]{QS}, which we mention only for comparison.
\end{remark}

\begin{corollary}\label{cor:degenerate}
Let $(R,\mu,\alpha)$ be a diagonal enhancement of a connected $R$ with effective data
$(Q,n,\varepsilon)$, and let $L$ be a link.
\begin{enumerate}[label=\textup{(\roman*)}]
\item If $n=0$, then $T_{R,\mu,\alpha}(L)=0$.
\item If $n=\pm1$, then $T_{R,\mu,\alpha}(L)=(\varepsilon n)^{|L|}$.
\item If $n=N\ge1$ and $\varepsilon=1$, then
$T_{R,\mu,\alpha}(L)=[N]_Q\,P_L(Q^N,Q-Q^{-1})$, the value at $q=Q$ of Turaev's invariant
$P_N(L)$ of \cite[Theorem~4.2.1 and Remark~4.2.6]{Tur}.
\item If $a=b$, then $Q^2=-1$ and $n=\ell$; the invariant vanishes if $\ell$ is even, and
$T_{R,\mu,\alpha}(L)=(\varepsilon Q^{\ell-1})^{|L|}$ with $\varepsilon Q^{\ell-1}=\pm1$ if
$\ell$ is odd.
\end{enumerate}
\end{corollary}

\begin{proof}
All four statements follow from Theorem~\ref{thm:connected-invariant}. (i) holds because
$[0]_Q=0$, and (iii) is the substitution $n=N$, $\varepsilon=1$. For (ii) we prove, more
generally, the identity
\begin{equation}\label{eq:rank-one}
 \varepsilon[n]_Q\,P_L(\varepsilon Q^n,Q-Q^{-1})=t^{|L|},\qquad t:=\varepsilon n,
\end{equation}
for all $Q$ with $Q^2\neq1$, $n=\pm1$ and $\varepsilon=\pm1$. Here $A=\varepsilon Q^n$
satisfies $A-A^{-1}=tz$ with $z=Q-Q^{-1}$, and $\varepsilon[n]_Q=t$. Since
$|L_+|=|L_-|=|L_0|\pm1$ and $t^2=1$, the invariant $L\mapsto t^{|L|}$ satisfies
$At^{|L_+|}-A^{-1}t^{|L_-|}=tz\,t^{|L_+|}=z\,t^{|L_0|}$, and on unlinks it agrees with the
left-hand side of~\eqref{eq:rank-one}, which is $t\,(tz/z)^{r-1}=t^r$ on $\bigcirc_r$. Hence Lemma~\ref{lem:skein-uniqueness} gives~\eqref{eq:rank-one}. For (iv), $Q^2=-a/b=-1$ and
$n=\ell$ by definition. Since $Q^{-1}=-Q$, we have
$[\ell]_Q=Q^{\ell-1}\bigl(1-(-1)^\ell\bigr)/2$, which vanishes for $\ell$ even. For $\ell$ odd
put $\varepsilon':=\varepsilon Q^{\ell-1}\in\{\pm1\}$; then $\varepsilon[\ell]_Q=\varepsilon'$ and
$\varepsilon Q^\ell=\varepsilon'Q$, and~\eqref{eq:rank-one} with $(Q,1,\varepsilon')$ gives the
claim.
\end{proof}

\begin{example}\label{ex:hecke-invariant}
For the standard Hecke operator with the enhancement $(R,\mu,Q^N)$ of
Example~\ref{ex:hecke-enhancement} we may take $d=1$. Then the effective data are $(Q,N,1)$ and
$T_{R,\mu,Q^N}(L)=[N]_Q\,P_L(Q^N,Q-Q^{-1})$. For $N=2$ this is, up to normalization and a
change of variable, the Jones polynomial \cite[Remark~4.3.6(iii)]{Tur}.
\end{example}

\subsection{State sums}
\label{subsec:state-sums}
Let $(R,\mu,\alpha)$ be a diagonal enhancement. By Proposition~\ref{prop:decoupling} each
$(R^{(s)},\mu^{(s)},\alpha)$ is a diagonal enhancement of the connected operator $R^{(s)}$.
Let $(Q_s,n_s,\varepsilon_s)$ be effective data for it, and put
\begin{equation}\label{eq:colour-data}
 A_s:=\varepsilon_sQ_s^{n_s},\qquad z_s:=Q_s-Q_s^{-1},\qquad \delta_s:=\varepsilon_s[n_s]_{Q_s},\qquad\omega_{st}:=\frac{\kappa_{st}}{\alpha^2}\ (s\neq t).
\end{equation}
Note that $\delta_s=(A_s-A_s^{-1})/z_s$. For a link $L$ and a colouring
$\chi\colon\pi_0(L)\to\{1,\dots,m\}$, let $L^\chi_s$ be the sublink of components of colour $s$.

The following lemma is the heart of the state sum. It makes precise the idea explained at the beginning of Section~\ref{subsec:colourings}. If $R$ is a braided direct sum, then every tensor factor keeps its colour as it moves along a strand, and therefore the braid acts on each colour separately, up to scalars coming from
the mixed crossings.

\begin{lemma}\label{lem:colour-words}
Assume that $R=R^{(1)}\bsum_\eta\dotsb\bsum_\eta R^{(m)}$, let $\mu$ be diagonal, and let
$u=\sigma_{k_1}^{\epsilon_1}\cdots\sigma_{k_q}^{\epsilon_q}$ be a word for $\xi\in B_r$. For
$c\in\{1,\dots,m\}^r$ put $V_c:=V_{\Lambda_{c_1}}\otimes\dots\otimes V_{\Lambda_{c_r}}$. Then
\[
 \tr\bigl(\rho_{R,r}(\xi)\mu^{\otimes r}\bigr)=\sum_c\ \prod_{s<t}\eta_{st}^{w_{st}(c)}
 \prod_{s\,:\,S_s(c)\neq\varnothing}
 \tr\Bigl(\rho_{R^{(s)},|S_s(c)|}\bigl(u_s(c)\bigr)\,(\mu^{(s)})^{\otimes|S_s(c)|}\Bigr),
\]
where $c$ runs over the words with $\pi_\xi\cdot c=c$, and $S_s(c):=\{i:c_i=s\}$, the words
$u_s(c)$ and the integers $w_{st}(c)$ are those of Section~\ref{subsec:colourings} and
Lemma~\ref{lem:deletion}.
\end{lemma}

\begin{proof}
\emph{Step 1. Braids permute the colour sectors.} We have $V^{\otimes r}=\bigoplus_cV_c$, the
sum over all $c\in\{1,\dots,m\}^r$. For $1\le k<r$ and every $c$, the operator
$R_k^{\pm1}$ maps $V_c$ onto $V_{\pi_{\sigma_k}\cdot c}$. Indeed, by Definition~\ref{def:braided-sum}, if
$c_k=c_{k+1}=s$ it acts on $V_c$ by $(R^{(s)})^{\pm1}$ on the factors $k,k+1$, and if
$c_k=s\neq t=c_{k+1}$ it acts as $\eta_{st}^{\pm1}$ times the exchange of these two factors.
Hence $\rho_{R,r}(\xi)=R_{k_1}^{\epsilon_1}\cdots R_{k_q}^{\epsilon_q}$ maps $V_c$ onto
$V_{\pi_\xi\cdot c}$, since $\pi_\xi=\pi_{\sigma_{k_1}}\cdots\pi_{\sigma_{k_q}}$.

\emph{Step 2. Only compatible colour words contribute.} Since $\mu^{\otimes r}$ preserves every
$V_c$, the operator
$\rho_{R,r}(\xi)\mu^{\otimes r}$ permutes the summands $V_c$, and only the summands with
$\pi_\xi\cdot c=c$ contribute to its trace, which gives
$\tr\bigl(\rho_{R,r}(\xi)\mu^{\otimes r}\bigr)=\sum_{\pi_\xi\cdot c=c}
\tr\bigl(\rho_{R,r}(\xi)\mu^{\otimes r}|_{V_c}\bigr)$.
Indeed, if $\pi_\xi\cdot c\neq c$, the operator maps $V_c$ into a different summand, and therefore the corresponding diagonal block of its matrix vanishes.

\emph{Step 3. Factorization on a compatible sector.} Fix $c$ with $\pi_\xi\cdot c=c$, let
$c^{(q)}=c,\dots,c^{(0)}=c$ be the words of colours at the levels of $D_u$
(Section~\ref{subsec:colourings}), and let $r_s:=|S_s(c)|$. The $p$-th letter maps
$V_{c^{(p)}}$ onto $V_{c^{(p-1)}}$, and
$\rho_{R,r}(\xi)|_{V_c}=R_{k_1}^{\epsilon_1}|_{V_{c^{(1)}}}\circ\dots\circ R_{k_q}^{\epsilon_q}|_{V_{c^{(q)}}}$.
For every colour word $c'$ with the same multiplicities as $c$ let
\[
 G_{c'}\colon V_{c'}\longrightarrow\bigotimes_{s\,:\,r_s>0}V_{\Lambda_s}^{\otimes r_s}
\]
be the isomorphism that groups the tensor factors by colour, keeping their order within each
colour. If the $p$-th crossing is monochromatic of colour $s$, then $c^{(p-1)}=c^{(p)}$, and
$G_{c^{(p-1)}}R_{k_p}^{\epsilon_p}G_{c^{(p)}}^{-1}$ is $\bigl(R^{(s)}_{j_p}\bigr)^{\epsilon_p}$ on the factor
$V_{\Lambda_s}^{\otimes r_s}$ and the identity on the other factors, where $j_p$ is the rank of
$k_p$ among the positions of colour $s$ at level $p$. If the $p$-th crossing is mixed, between
colours $s$ and $t$, then $G_{c^{(p-1)}}R_{k_p}^{\epsilon_p}G_{c^{(p)}}^{-1}=\eta_{st}^{\epsilon_p}\,\id$,
because exchanging two adjacent factors of different colours does not change the order of the
factors within each colour. Composing these identities for $p=q,\dots,1$, the intermediate maps
$G_{c^{(p)}}$ cancel. The monochromatic letters of colour $s$, in their order in $u$ and with
the indices $j_p$, form the word $u_s(c)$, and the scalars of the mixed letters multiply to
$\prod_{s<t}\eta_{st}^{w_{st}(c)}$. Hence
\[
 G_c\,\rho_{R,r}(\xi)|_{V_c}\,G_c^{-1}=\prod_{s<t}\eta_{st}^{w_{st}(c)}
 \bigotimes_{s\,:\,r_s>0}\rho_{R^{(s)},r_s}\bigl(u_s(c)\bigr),\qquad
 G_c\,\mu^{\otimes r}|_{V_c}\,G_c^{-1}=\bigotimes_{s\,:\,r_s>0}(\mu^{(s)})^{\otimes r_s}.
\]

\emph{Step 4. The trace.} Conjugating by $G_c$ does not change the trace, and the trace is
multiplicative under tensor products. Hence the contribution of $c$ is
$\prod_{s<t}\eta_{st}^{w_{st}(c)}\prod_{s\,:\,r_s>0}\tr\bigl(\rho_{R^{(s)},r_s}(u_s(c))
(\mu^{(s)})^{\otimes r_s}\bigr)$, which is the term of $c$ in the claimed formula.
\end{proof}

For the word $u=\sigma_1\sigma_1\sigma_2$ and the colour word $c=(1,2,2)$ of
Example~\ref{ex:coloured-braid} (Figure~\ref{fig:coloured-braid}), we have $w_{12}(c)=2$, the
word $u_1(c)$ is empty and $u_2(c)=\sigma_1\in B_2$. Hence the term of $c$ in
Lemma~\ref{lem:colour-words} is
$\eta_{12}^2\,\tr(\mu^{(1)})\,\tr\bigl(R^{(2)}(\mu^{(2)}\otimes\mu^{(2)})\bigr)$.

\begin{theorem}[Theorem~\ref{mainthm:C}]\label{thm:state-sum}
Let $(R,\mu,\alpha)$ be a diagonal enhancement. For every oriented link $L$,
\begin{align*}
 T_{R,\mu,\alpha}(L)
 &=\sum_{\chi}\ \prod_{s<t}\omega_{st}^{\lk(L^\chi_s,L^\chi_t)}
 \prod_{s\,:\,L^\chi_s\neq\varnothing}T_{R^{(s)},\mu^{(s)},\alpha}(L^\chi_s)\\
 &=\sum_{\chi}\ \prod_{s<t}\omega_{st}^{\lk(L^\chi_s,L^\chi_t)}
 \prod_{s\,:\,L^\chi_s\neq\varnothing}\delta_s\,P_{L^\chi_s}(A_s,z_s),
\end{align*}
where $\chi$ runs over all maps $\pi_0(L)\to\{1,\dots,m\}$.
\end{theorem}

\begin{proof}
\emph{Step 1. Reduction to a braided direct sum.} Choose $\eta_{st}\in\C^\times$ with $\eta_{st}=\eta_{ts}$ and $\eta_{st}^2=\kappa_{st}$. By Proposition~\ref{thm:canonical-decomposition} and
Proposition~\ref{prop:px-invariance}, we may replace $R$ by $R^{(1)}\bsum_\eta\dotsb\bsum_\eta R^{(m)}$ without changing $\mu$, $\alpha$ or
the invariant. The X-transformation used there fixes the restrictions $R^{(s)}$, and the
scalars $\kappa_{st}$ are X-invariant; hence the right-hand sides do not change either.

\emph{Step 2. Colourings as colour words.} Write $L=\widehat\xi$ with $\xi\in B_r$ and choose
a word $u$ for $\xi$. The colourings $\chi$ of $L$ correspond bijectively to the compatible
colour words $c$, by $c_i:=\chi(\text{component through the bottom position }i)$. Then $S_s(c)\neq\varnothing$ if and
only if $L^\chi_s\neq\varnothing$, and by Lemma~\ref{lem:deletion} the closure of $u_s(c)$ is
$L^\chi_s$, $w_{st}(c)=2\lk(L^\chi_s,L^\chi_t)$ and
$w(\xi)=\sum_{s\,:\,S_s(c)\neq\varnothing}w(u_s(c))+\sum_{s<t}w_{st}(c)$.

\emph{Step 3. The state sum.} Put $r_s:=|S_s(c)|$. We multiply the identity of
Lemma~\ref{lem:colour-words} by $\alpha^{-w(\xi)}$ and distribute the factor
$\alpha^{-w(\xi)}=\prod_{s\,:\,r_s>0}\alpha^{-w(u_s(c))}\prod_{s<t}\alpha^{-w_{st}(c)}$ over the term of
$c$. Using $w_{st}(c)=2\lk(L^\chi_s,L^\chi_t)$, the term of $c$ becomes
\begin{align*}
 &\prod_{s<t}\Bigl(\frac{\eta_{st}}{\alpha}\Bigr)^{2\lk(L^\chi_s,L^\chi_t)}
 \prod_{s\,:\,L^\chi_s\neq\varnothing}\alpha^{-w(u_s(c))}
 \tr\Bigl(\rho_{R^{(s)},r_s}\bigl(u_s(c)\bigr)(\mu^{(s)})^{\otimes r_s}\Bigr)\\
 &\qquad=\prod_{s<t}\omega_{st}^{\lk(L^\chi_s,L^\chi_t)}
 \prod_{s\,:\,L^\chi_s\neq\varnothing}T_{R^{(s)},\mu^{(s)},\alpha}(L^\chi_s),
\end{align*}
by~\eqref{eq:invariant} applied to the diagonal enhancement $(R^{(s)},\mu^{(s)},\alpha)$ and
the braid represented by $u_s(c)$, whose closure is $L^\chi_s$; here we used
$\eta_{st}^2/\alpha^2=\kappa_{st}/\alpha^2=\omega_{st}$. Summing over $c$, that is, over $\chi$,
gives the first equality. The second follows from Theorem~\ref{thm:connected-invariant} applied
to each $(R^{(s)},\mu^{(s)},\alpha)$, and Remark~\ref{rem:choices} shows that the result does
not depend on the choice of the effective data.
\end{proof}

\begin{example}\label{ex:hopf-state-sum}
Let $H_+$ be the positive Hopf link, the closure of $\sigma_1^2\in B_2$; its components are
unknots with linking number $1$. For a colouring $\chi$, either both components have the same
colour $s$, which contributes $T_{R^{(s)},\mu^{(s)},\alpha}(H_+)=\delta_sP_{H_+}(A_s,z_s)$, or
they have different colours $s\neq t$, which contributes $\omega_{st}\delta_s\delta_t$. Hence
\[
 T_{R,\mu,\alpha}(H_+)=\sum_{s=1}^m\delta_s\,P_{H_+}(A_s,z_s)+
 \sum_{s\neq t}\omega_{st}\,\delta_s\delta_t .
\]
It is instructive to see the mixed terms directly in the trace. For $R$ a braided direct sum
and $c=(s,t)$ with $s\neq t$, the operator $R^2=\rho_{R,2}(\sigma_1^2)$ acts on
$V_{\Lambda_s}\otimes V_{\Lambda_t}$ as $\eta_{st}^2\,\id$, because $R$ exchanges the two factors
and multiplies by $\eta_{st}$ each time. Hence this sector contributes
$\alpha^{-2}\eta_{st}^2\tr(\mu^{(s)})\tr(\mu^{(t)})=\omega_{st}\delta_s\delta_t$ to
$T_{R,\mu,\alpha}(H_+)=\alpha^{-2}\tr\bigl(R^2(\mu\otimes\mu)\bigr)$.
\end{example}

\begin{corollary}\label{cor:state-sum-consequences}
Let $(R,\mu,\alpha)$ be a diagonal enhancement.
\begin{enumerate}[label=\textup{(\roman*)}]
\item For every knot $K$, $T_{R,\mu,\alpha}(K)=\sum_{s=1}^m\delta_s\,P_K(A_s,z_s)$.
\item If all pairwise linking numbers of the components of $L$ vanish, then
$T_{R,\mu,\alpha}(L)$ does not depend on the $\kappa_{st}$.
\item Let $f\colon\pi_0(L)\to\pi_0(L')$ be a bijection that preserves the pairwise linking
numbers and satisfies $P_{f(J)}=P_J$ for every sublink $J$ of $L$, where $f(J)$ is the
corresponding sublink of $L'$. Then $T_{R,\mu,\alpha}(L)=T_{R,\mu,\alpha}(L')$.
\end{enumerate}
\end{corollary}

\begin{proof}
(i) A knot $K$ has one component, hence its colourings are the $m$ constant maps. For the
colouring with value $s$, the only nonempty sublink is $K$ itself, of colour $s$, and every
linking number in Theorem~\ref{thm:state-sum} involves an empty sublink. Hence this colouring
contributes $\delta_sP_K(A_s,z_s)$.

(ii) By Section~\ref{subsec:braids}, $\lk(L^\chi_s,L^\chi_t)$ is the sum of the linking numbers
of the pairs formed by a component of colour $s$ and a component of colour $t$. If all pairwise
linking numbers vanish, every factor $\omega_{st}^{\lk(L^\chi_s,L^\chi_t)}$ is equal to $1$, and
the remaining factors do not involve the $\kappa_{st}$.

(iii) The map $\chi'\mapsto\chi'\circ f$ is a bijection from the colourings of $L'$ onto the
colourings of $L$. Let $\chi=\chi'\circ f$. Then $f(L^\chi_s)=(L')^{\chi'}_s$ for every $s$, hence
$P_{L^\chi_s}=P_{(L')^{\chi'}_s}$ by hypothesis, and
$\lk(L^\chi_s,L^\chi_t)=\lk\bigl((L')^{\chi'}_s,(L')^{\chi'}_t\bigr)$ by the description of these
numbers in~(ii), because $f$ preserves the pairwise linking numbers. Therefore the terms of
$\chi$ and of $\chi'$ in Theorem~\ref{thm:state-sum} coincide.
\end{proof}

\begin{example}\label{ex:flip-invariant}
Consider the scaled flip with the enhancement $(\lambda\tau,\id,\lambda)$ on $V$ of dimension
$N$. Each component is a singleton with $n_s=1$ and $\delta_s=1$, hence $T_{R^{(s)},\mu^{(s)},\lambda}(J)=1$ for every link $J$ by Corollary~\ref{cor:degenerate}(ii),
and $\omega_{st}=\lambda^2/\lambda^2=1$. Every colouring contributes $1$, hence $T_{\lambda\tau,\id,\lambda}(L)=N^{|L|}$, and $T_{\lambda\tau,-\id,-\lambda}(L)=(-N)^{|L|}$.
For $N\ge2$ the invariant $N^{|L|}$ is not of the form $\delta P_L(A,z)$ with
$A,z\in\C^\times$. Indeed, the unknot and the two-component unlink would force $\delta=N$ and
$A-A^{-1}=Nz$. The positive Hopf link $H_+$, the two-component unlink and the unknot form a
skein triple $(L_+,L_-,L_0)$, and the skein relation would give $(A-A^{-1})N^2=zN$, that is, $N^3z=Nz$, which is impossible because $z\neq0$ and $N\ge2$.
\end{example}

\emph{The gluing of Majid and Markl.} Majid and Markl~\cite{MM} introduced a gluing operation
$\oplus_q$ for Yang--Baxter operators of Hecke type, together with a gluing of the associated
link invariants. The mixed blocks of $\oplus_q$ are of Hecke type. Hence, for charge-conserving operators with the same spectral parameters, $\oplus_q$ concatenates the scalar classes into a
single connected component, and the glued invariant is a state sum over colourings of the edges
of a diagram in which colours may change at crossings \cite[Section~6]{MM}. The braided direct
sums of Definition~\ref{def:braided-sum} are of a different nature. Their summands may have
unrelated spectral parameters, the colours in Theorem~\ref{thm:state-sum} are constant along the
components of the link, and the resulting invariants are in general not specializations of the
HOMFLY--PT polynomial (Example~\ref{ex:flip-invariant}).

\subsection{Realization}
\label{subsec:realization}
We now show that every state sum of the form~\eqref{eq:state-sum-intro} occurs. By
Theorem~\ref{thm:normal-form}(b) and Proposition~\ref{prop:decoupling}, it suffices to realize
each colour separately by an operator of the form~\eqref{eq:normal-form} with $\alpha=1$, and
then to glue these operators by suitable multiples of the flip. The case $Q^2=-1$, where the
two spectral parameters coincide, needs a separate argument.

\begin{lemma}\label{lem:connected-realization}
Let $Q\in\C^\times$ with $Q^2\notin\{1,-1\}$, $n\in\Z$ and $\varepsilon\in\{\pm1\}$. There are
an indexed line decomposition $(V,\cL)$, an operator $R$ of the form~\eqref{eq:normal-form}
whose scalar classes are singletons, and a diagonal enhancement $(R,\mu,1)$ with effective data
$(Q,n,\varepsilon)$.
\end{lemma}

\begin{proof}
Choose $\ell\ge2$ with $\ell\ge|n|$ and $\ell\equiv n\pmod2$, and put $p:=(\ell+n)/2$, which satisfies $0\le p\le\ell$. Let $d:=\varepsilon Q^{-n}$, $a:=dQ$ and $b:=-dQ^{-1}$; then
$a+b=d(Q-Q^{-1})\neq0$, and $a\neq b$ because $Q^2\neq-1$. Theorem~\ref{thm:normal-form}(a),
applied to $I=\{1,\dots,\ell\}$ with the classes $\{1\}\prec\dots\prec\{\ell\}$, $v_r=a$ for
$r\le p$ and $v_r=b$ for $r>p$, gives a connected $R$ with spectral parameters $a,b$ and
signed effective rank $p-(\ell-p)=n$. Since $-ab=d^2$ and $-a/b=Q^2$, we get
$\cE(R)=d^2Q^{2n}=1$, and Theorem~\ref{thm:connected-existence} gives a diagonal enhancement
$(R,\mu,1)$. Computed with the ordering $(a,b)$ and the square root $d$ of $-ab$, its
effective data are $(a/d,\,n,\,1/(dQ^n))=(Q,n,\varepsilon)$.
\end{proof}

\begin{lemma}\label{lem:repeated-realization}
Let $Q\in\C^\times$ with $Q^2=-1$, $n\in\Z$ and $\varepsilon\in\{\pm1\}$. There is a diagonal
enhancement $(R,\mu,1)$ of an operator $R$ of the form~\eqref{eq:normal-form} with
$T_{R,\mu,1}(L)=\varepsilon[n]_Q\,P_L(\varepsilon Q^n,Q-Q^{-1})$ for every link $L$.
\end{lemma}

\begin{proof}
Since $Q^{-1}=-Q$, we have $[n]_Q=Q^{n-1}\bigl(1-(-1)^n\bigr)/2$. If $n$ is even, the
right-hand side vanishes; take $I=\{1,2\}$ and the operator~\eqref{eq:normal-form} with classes
$\{1\}\prec\{2\}$ and $a=b=v_1=v_2=Q$, which is connected by
Theorem~\ref{thm:normal-form}(a). Its signed effective rank is $2$ and
$\cE(R)=-Q^2=1$. Hence it has a diagonal enhancement $(R,\mu,1)$, whose invariant vanishes by
Corollary~\ref{cor:degenerate}(iv). If $n$ is odd, put $t:=\varepsilon Q^{n-1}\in\{\pm1\}$;
then $\varepsilon[n]_Q=t$ and $\varepsilon Q^n=tQ$, and by~\eqref{eq:rank-one} the right-hand side equals $t^{|L|}$. Take $\dim V=1$ and $R=t\,\id$, which is the
operator~\eqref{eq:normal-form} with one class and $a=b=v_1=t$, for which $\cE(R)=t^2=1$, and
take $\mu=t$; then~\eqref{eq:enhancement} holds with $\alpha=1$, and
$T_{R,\mu,1}(\widehat\xi)=t^{w(\xi)}t^{r}=t^{|\widehat\xi|}$ by~\eqref{eq:parity}.
\end{proof}

\begin{theorem}[Theorem~\ref{mainthm:D}]\label{thm:realization}
Let $m\ge1$; for $1\le s\le m$ let $Q_s\in\C^\times$ with $Q_s^2\neq1$, $n_s\in\Z$ and
$\varepsilon_s\in\{\pm1\}$, and for $s\neq t$ let $\omega_{st}\in\C^\times$ with $\omega_{st}=\omega_{ts}$. Define
$A_s,z_s,\delta_s$ as in~\eqref{eq:colour-data}. Then there is a diagonal enhancement
$(R,\mu,1)$ of a charge-conserving Yang--Baxter operator such that, for every link $L$,
\[
 T_{R,\mu,1}(L)=\sum_{\chi}\ \prod_{s<t}\omega_{st}^{\lk(L^\chi_s,L^\chi_t)}
 \prod_{s\,:\,L^\chi_s\neq\varnothing}\delta_s\,P_{L^\chi_s}(A_s,z_s).
\]
Consequently, the link invariants of diagonal enhancements of charge-conserving Yang--Baxter
operators are exactly the state sums of this form.
\end{theorem}

\begin{proof}
\emph{Step 1. One colour at a time.} For each $s$, Lemma~\ref{lem:connected-realization} (if
$Q_s^2\neq-1$) together with Theorem~\ref{thm:connected-invariant}, or
Lemma~\ref{lem:repeated-realization} (if $Q_s^2=-1$), gives an indexed line decomposition
$(V_s,\cL_s)$ with index set $I_s$, an operator $R^{(s)}$ of the form~\eqref{eq:normal-form} on
$(V_s,\cL_s)$ and a diagonal enhancement $(R^{(s)},\mu^{(s)},1)$ with $T_{R^{(s)},\mu^{(s)},1}(J)=\delta_sP_J(A_s,z_s)$ for
every link $J$.

\emph{Step 2. Gluing.} Let $V:=\bigoplus_sV_s$, with the union $\cL$ of the $\cL_s$, indexed by $I:=\bigsqcup_sI_s$; choose $\eta_{st}\in\C^\times$ with $\eta_{st}=\eta_{ts}$ and $\eta_{st}^2=\omega_{st}$, and let $R:=R^{(1)}\bsum_\eta\dotsb\bsum_\eta R^{(m)}$ and $\mu:=\bigoplus_s\mu^{(s)}$. By
Theorem~\ref{thm:normal-form}(b), $R\in\YBCC(V,\cL)$, its connected components are the $I_s$,
its connected restrictions are the $R^{(s)}$ and $\kappa_{st}=\omega_{st}$. By
Proposition~\ref{prop:decoupling}, $(R,\mu,1)$ is a diagonal enhancement, since each
$(R^{(s)},\mu^{(s)},1)$ is one.

\emph{Step 3. The invariant.} With $\alpha=1$ we have $\kappa_{st}/\alpha^2=\omega_{st}$. Hence the first formula of Theorem~\ref{thm:state-sum}, combined with Step~1, gives the claimed formula.

\emph{Step 4. The converse.} By Theorem~\ref{thm:state-sum} and Lemma~\ref{lem:effective-data},
the invariant of any diagonal enhancement is a state sum of this form, with $Q_s^2\neq1$,
$\varepsilon_s=\pm1$ and $\omega_{st}=\omega_{ts}\in\C^\times$.
\end{proof}

\section{Examples and comparisons}
\label{sec:examples}

In this section we compare our results with the literature. In dimension two they recover the
results of Aizawa, Harada, Kawaguchi and Otsuki~\cite{AHKO} for the charge-conserving solutions
(Section~\ref{subsec:dim-two}), and they explain the invariants of the dressed solutions computed
in \cite[Section~5]{AHKO} (Section~\ref{subsec:dressings}). For the operators through which the
Yokonuma--Hecke algebras act on tensor space~\cite{LPdA}, Theorem~\ref{thm:state-sum} and
Lickorish's formula \cite[Appendix~B]{CJKL} identify the invariant with a specialization of the
three-variable invariant $\Theta$ of Chlouveraki, Juyumaya, Karvounis and
Lambropoulou~\cite{CJKL} (Section~\ref{subsec:YH}).

\subsection{Dimension two}
\label{subsec:dim-two}
Let $I=\{1,2\}$. By Proposition~\ref{prop:two-index}, which for $I=\{1,2\}$ is the
classification of the two-dimensional case in \cite[Proposition~3.22]{MR}, up to PX-equivalence
an operator $R\in\YBCC(V,\cL)$ is of one of the following types, namely $R=v\,\id$ (type $0$, the edge
$\{1,2\}$ is scalar); the normal form with classes $\{1\}\prec\{2\}$ and class scalars
$(a,a)$ (type $\mathrm f$) or $(a,b)$ with $a\neq b$ (type $\mathrm a$); or the braided direct
sum of the one-dimensional operators given by the vertex scalars $v_1$ and $v_2$, with coupling $\eta$ (type $/$). (Class
scalars $(b,b)$ and $(b,a)$ reduce to types $\mathrm f$ and $\mathrm a$ by exchanging the names
of $a$ and $b$.) Theorems~\ref{thm:connected-existence}, \ref{thm:global-existence},
\ref{thm:connected-invariant} and~\ref{thm:state-sum} give the following table.
\begin{center}
\renewcommand{\arraystretch}{1.3}
\begin{tabular}{c|c|c|c|c}
type & $(\ell,n)$ & $\alpha^2$ & $(\mu_1,\mu_2)$ & $T_{R,\mu,\alpha}(L)$\\ \hline
$0$ & $(1,1)$ & $v^2$ & $\mu_1+\mu_2=v/\alpha$ & $(v/\alpha)^{|L|}$\\
$\mathrm f$ & $(2,2)$ & $-a^3/b$ & $\bigl(a/\alpha,\,-a^2/(b\alpha)\bigr)$ &
$\varepsilon[2]_Q\,P_L(\varepsilon Q^2,Q-Q^{-1})$\\
$\mathrm a$ & $(2,0)$ & $-ab$ & $\bigl(a/\alpha,\,-a/\alpha\bigr)$ & $0$\\
$/$ & $(1,1),(1,1)$ & $v_1^2=v_2^2$ & $\bigl(v_1/\alpha,\,v_2/\alpha\bigr)$ &
$\sum_\chi\omega^{\lk(L^\chi_1,L^\chi_2)}\delta_1^{|L^\chi_1|}\delta_2^{|L^\chi_2|}$
\end{tabular}
\end{center}
Here $(Q,2,\varepsilon)$ are effective data in type $\mathrm f$, and in type $/$ we have
$\delta_s=v_s/\alpha\in\{\pm1\}$ and $\omega=\eta^2/\alpha^2$. For $a=Q$ and $b=-Q^{-1}$, type
$\mathrm f$ is the standard Hecke operator with $N=2$. Type $/$ admits a diagonal enhancement
only if $v_1=\pm v_2$, whatever $\eta$ is; for the Hopf links $H_\pm$, with linking number
$\pm1$, it gives $T_{R,\mu,\alpha}(H_\pm)=2+2\,\delta_1\delta_2\,\omega^{\pm1}$.

In dimension two, an invertible diagonal $\mu$ is generic in the sense of the introduction,
that is, $\mu_i\mu_j\neq\mu_k\mu_l$ whenever $\{i,j\}\neq\{k,l\}$, if and only if
$\mu_1\neq\pm\mu_2$. In the table this happens only in type $0$, for suitable $\mu$, and in
type $\mathrm f$ with $a\neq b$, where $\mu_2/\mu_1=-a/b$. By
the discussion leading to~\eqref{eq:cc-intro} and Lemmas~\ref{lem:conjugation}
and~\ref{lem:turaev-relation},
the invariant of every enhanced Yang--Baxter operator on $\C^2$ whose enhancing endomorphism is
diagonalizable with nonzero eigenvalues $\mu_1\neq\pm\mu_2$ is therefore either $(\pm1)^{|L|}$ or a
specialization $\varepsilon[2]_Q\,P_L(\varepsilon Q^2,Q-Q^{-1})$ of the Jones polynomial.

All solutions of the Yang--Baxter equation in dimension two were found by
Hietarinta~\cite{Hie92,Hie93}, and Aizawa, Harada, Kawaguchi and Otsuki~\cite{AHKO} determined
the enhancements of the nonsingular solutions in his list that depend on at least one
parameter. In their notation \cite[Section~3]{AHKO}, with the basis $(e_+,e_-)$, the
charge-conserving solutions in this list are $R_{3.1}$ (type $/$ with $v_1=1$, $v_2=s$ and
$\kappa_{12}=pq$), $R_{2.1}$ (type $\mathrm f$ with classes $\{+\}\prec\{-\}$, $a=1$ and
$b=-pq$) and $R_{2.2}$ (type $\mathrm a$ with $a=1$ and $b=-pq$). For invertible $\mu$, the
enhancements of these three solutions listed in \cite[Table~1]{AHKO} are exactly those of the
table above; their enhancing endomorphism is $\beta$ times ours
(Lemma~\ref{lem:turaev-relation}). For instance, for $R_{2.1}$ one finds
$\alpha=\pm1/\sqrt{pq}$ and $\mu=\pm(\sqrt{pq},1/\sqrt{pq})$, with the specialization of the
Jones polynomial for $Q=1/\sqrt{pq}$, and for $R_{2.2}$ one finds $\alpha=\pm\sqrt{pq}$ and
$\mu=\pm(1,-1)/\sqrt{pq}$, with vanishing invariant. The other enhancements of these solutions in \cite[Table~1]{AHKO} have non-invertible $\mu$. Those with diagonal $\mu$ are the ones of Remark~\ref{rem:singular-mu} with $J=\{i\}$ a single index, namely $\alpha=\pm v_i$, where $v_i$ is the vertex scalar, and $\mu_i=\beta\alpha/v_i$; their invariant is $(\alpha/v_i)^{|L|}$. The remaining ones occur for $R_{2.1}$ with $q=1$, and their $\mu$ is not diagonal.

\subsection{Dressings}
\label{subsec:dressings}
Aizawa et al.\ also considered \emph{diagonal dressings} \cite[Section~5]{AHKO}. Given a
solution $\widetilde R$ on $V_J$, $J\subset I$, the dressed operator $R$ on $V$ agrees with
$\widetilde R$ on $V_J^{\otimes2}$ and satisfies $R(e_i\otimes e_j)=s_{ji}\,e_j\otimes e_i$
whenever $i\notin J$ or $j\notin J$ \cite[(5.1)]{AHKO}. If $\widetilde R$ is charge-conserving, then $R$ is also charge-conserving; every $k\notin J$ is then a one-element connected component of $R$, and by
Proposition~\ref{thm:canonical-decomposition} $R$ is X-equivalent to the braided direct sum of the
connected restrictions of $\widetilde R$ and the one-dimensional operators $s_{kk}$,
$k\notin J$. Hence Theorem~\ref{thm:state-sum} describes the invariants of all such dressings.
On a knot, the dressed invariant is the invariant of $\widetilde R$ plus the constant
$\sum_{k\notin J}\mu_k$, and the parameters $s_{ij}$ with $i\neq j$ enter only through the
couplings $\omega_{st}$, that is, through linking numbers; this is what is observed in
\cite[Section~5.2]{AHKO}. For example, for the dressing of $R_{2.1}$ in \cite[(5.28)]{AHKO},
with $J=\{1,3\}$, $\mu=(\sqrt{pq},1,1/\sqrt{pq})$ and $\alpha=1/\sqrt{pq}$, there are two
colours, $\{1,3\}$ and $\{2\}$, with $\kappa_{12}=s_{12}s_{21}$ and $\omega_{12}=\kappa_{12}/\alpha^2=s_{12}s_{21}\,pq$. In \cite[Tables~2 and~3]{AHKO} the product $s_{12}s_{21}$ is denoted $s$ and $pq$ is denoted $t$, hence $\omega_{12}$ is the product $st$ that appears there. The colour $\{2\}$ carries the one-dimensional
operator $s_{22}=\alpha$ with $\mu_2=1$, whose invariant is identically $1$
(Corollary~\ref{cor:degenerate}(ii)). Hence Theorem~\ref{thm:state-sum} gives
\[
 T_{R,\mu,\alpha}(L)=\sum_{\chi\colon\pi_0(L)\to\{1,2\}}(s_{12}s_{21}\,pq)^{\lk(L^\chi_1,L^\chi_2)}\,
 T_{R_{2.1},\widetilde\mu,\alpha}(L^\chi_1),\qquad T_{R_{2.1},\widetilde\mu,\alpha}(\varnothing)=1.
\]
For the dressings of $R_{2.2}$ in \cite[(5.29) and (5.30)]{AHKO}, the invariant of $R_{2.2}$
vanishes on every nonempty link (Corollary~\ref{cor:degenerate}(i)), and therefore only the colourings
that avoid the colour of $R_{2.2}$ contribute. We have checked by direct computation that the
values in \cite[Tables~2--4]{AHKO} (the first and the last normalized by the value on the
unknot) agree with Theorem~\ref{thm:state-sum}, except for an apparent misprint in the entry
for $6^2_3$ in \cite[Table~3]{AHKO}, where the coefficient of $t^7$ should be $1$, as also
follows from the Jones polynomial given in the same row.

\subsection{Yokonuma--Hecke algebras}
\label{subsec:YH}
We first recall the objects involved, following \cite{CJKL}; in this subsection the letters
$q$, $\lambda$ and $E$ denote the variables of \cite{CJKL}. Let $m,r\ge1$ and
$q\in\C^\times$. The modular framed braid group $(\Z/m\Z)^r\rtimes B_r$ is generated by
$\sigma_1,\dots,\sigma_{r-1}$ and by commuting elements $t_1,\dots,t_r$, subject to the braid
relations, $t_j^m=1$ and $t_j\sigma_k=\sigma_kt_{\pi_{\sigma_k}(j)}$. The
\emph{Yokonuma--Hecke algebra} $Y_{m,r}(q)$ is the quotient of its group algebra by the
relations $g_k^2=1+(q-q^{-1})f_kg_k$, where $g_k$ is the image of $\sigma_k$ and $f_k:=\frac1m\sum_{u=0}^{m-1}t_k^ut_{k+1}^{-u}$, which is an idempotent \cite[Section~1]{CJKL}; for $m=1$ it is
the Iwahori--Hecke algebra. Markov traces on these algebras were constructed by Juyumaya \cite{Juy04}.
When the parameters of such a trace satisfy the E-condition of \cite{JL13}, the trace can be rescaled to give invariants of framed links, and their restrictions to classical links give invariants
$\Theta_m$ \cite{CJKL}. Chlouveraki, Juyumaya, Karvounis and Lambropoulou showed that the
$\Theta_m$ are specializations of a three-variable invariant $\Theta(q,\lambda,E)$ of oriented
links, namely $\Theta_m$ is $\Theta$ at $E=1/m$ \cite[Section~8]{CJKL}. The invariant
$\Theta$ is characterized by two rules \cite[Theorem~8.1]{CJKL}. First, $\Theta(L)=E^{1-k}P(L)$
if $L$ is a disjoint union of $k$ knots, where $P$ is the HOMFLY--PT polynomial normalized by
$\lambda^{-1/2}P(L_+)-\lambda^{1/2}P(L_-)=(q-q^{-1})P(L_0)$ and $P(\bigcirc)=1$. Second, the
same skein relation $\lambda^{-1/2}\Theta(L_+)-\lambda^{1/2}\Theta(L_-)=(q-q^{-1})\Theta(L_0)$ holds
for skein triples whose two crossing strands belong to different components.

When all colours carry the same data, the state sum of Theorem~\ref{thm:state-sum} can be
grouped according to partitions of the set of components, and this produces $\Theta$.

\begin{proposition}\label{prop:YH}
Let $(R,\mu,\alpha)$ be a diagonal enhancement of an operator with $m$ connected components,
and suppose that $(A_s,z_s,\delta_s)=(A,z,\delta)$ for all $s$ and $\omega_{st}=\omega$ for all
$s\neq t$. For a link $L$ and a partition $\pi$ of $\pi_0(L)$, let $\nu(\pi)$ be the sum of the
linking numbers of the pairs of components that lie in different blocks of $\pi$, and let $P_\pi(L):=\prod_{B\in\pi}P_{L_B}(A,z)$, where $L_B$ is the sublink formed by the components in
$B$. Then
\begin{equation}\label{eq:grouped}
 T_{R,\mu,\alpha}(L)=\sum_{\pi}m(m-1)\cdots(m-|\pi|+1)\,\delta^{|\pi|}\,\omega^{\nu(\pi)}\,
 P_\pi(L),
\end{equation}
the sum running over all partitions $\pi$ of $\pi_0(L)$. If moreover $\omega=A^{-2}$, then
$T_{R,\mu,\alpha}(L)=m\delta\,\Theta(L)$, where $\Theta$ is the invariant of
\cite[Theorem~8.1]{CJKL} with its variables specialized to $\lambda^{-1/2}=A$,
$q-q^{-1}=z$ and $E=1/m$.
\end{proposition}

\begin{proof}
Group the colourings $\chi$ in Theorem~\ref{thm:state-sum} according to the partition $\pi$
of $\pi_0(L)$ into the nonempty fibres of $\chi$. A partition with $k$ blocks arises from
exactly $m(m-1)\cdots(m-k+1)$ colourings, namely the injective maps from the set of blocks to
$\{1,\dots,m\}$. Each of them contributes $\omega^{\nu(\pi)}\delta^kP_\pi(L)$, because two
components have different colours if and only if they lie in different blocks. This
proves~\eqref{eq:grouped}.

Lickorish's formula \cite[Theorem~B.1]{CJKL} states that
\[
 \Theta(L)=\sum_{k\ge1}\Bigl(\frac{\lambda^{-1/2}-\lambda^{1/2}}{q-q^{-1}}\Bigr)^{k-1}E_k
 \sum_{|\pi|=k}\lambda^{\nu(\pi)}P(\pi L),\qquad E_k=(E^{-1}-1)(E^{-1}-2)\cdots(E^{-1}-k+1),
\]
where $P(\pi L)$ is the product of the polynomials $P(L_B)$ of the sublinks $L_B$, $B\in\pi$,
with $P$ normalized as above. With $\lambda^{-1/2}=A$ and $q-q^{-1}=z$ this is $P_\pi(L)$, the factor
$(\lambda^{-1/2}-\lambda^{1/2})/(q-q^{-1})=(A-A^{-1})/z$ equals $\delta$, and for $E=1/m$ we
have $mE_k=m(m-1)\cdots(m-k+1)$. Comparing with~\eqref{eq:grouped} for $\omega=A^{-2}=\lambda$
gives $T_{R,\mu,\alpha}(L)=m\delta\,\Theta(L)$.
\end{proof}

\begin{example}\label{ex:YH}
Let $N\ge2$ and $Q\in\C^\times$ with $Q^2\neq1$, and let $R$ be the braided direct sum of $m$
copies of the standard Hecke operator on $\C^N$, with $\eta_{st}=1$ for all $s\neq t$. Up to
PX-equivalence, which does not change the invariants (Proposition~\ref{prop:px-invariance}),
this is the operator $\sigma$ of Lacabanne and Poulain d'Andecy \cite[Section~3.2]{LPdA} for
$V_1=\dots=V_m=\C^N$, through which the Yokonuma--Hecke algebra acts on $V^{\otimes r}$; for
generic $Q$ this action maps onto the centralizer of
$U_Q(\mathfrak{gl}_N\oplus\dots\oplus\mathfrak{gl}_N)$ \cite[Theorem~5.2]{LPdA}. By Example~\ref{ex:hecke-enhancement} and
Proposition~\ref{prop:decoupling}, $(R,\mu,Q^N)$ is a diagonal enhancement for $\mu$ the direct
sum of $m$ copies of $\operatorname{diag}(Q^{1-N},Q^{3-N},\dots,Q^{N-1})$. Every colour carries
the data $(A,z,\delta)=(Q^N,Q-Q^{-1},[N]_Q)$, and $\omega_{st}=Q^{-2N}=A^{-2}$.
Proposition~\ref{prop:YH} gives
\[
 T_{R,\mu,Q^N}(L)=m[N]_Q\,\Theta(L)\qquad\text{at }\lambda=Q^{-2N},\ q=Q,\ E=1/m .
\]
For $E=1/m$, the invariant $\Theta$ is the invariant $\Theta_m$ obtained from the Markov
traces on the Yokonuma--Hecke algebras \cite[Section~8]{CJKL}.
\end{example}

Formulas of the same shape as~\eqref{eq:grouped}, in terms of the HOMFLY--PT polynomials of
sublinks and linking numbers, were obtained for the invariants of classical links defined by
Markov traces on the Yokonuma--Hecke algebras~\cite{Juy04,JL13,ChL13}, namely Lickorish's
closed formula \cite[Appendix~B]{CJKL}, the work of Jacon and Poulain
d'Andecy~\cite{JPdA16}, and the description of all these invariants by Poulain d'Andecy and
Wagner~\cite{PdAW}; see also~\cite{KL}. Poulain d'Andecy and Wagner compute the image of a
braid in terms of the braids obtained by deleting strands \cite[Theorem~4.1]{PdAW}, an
algebraic counterpart of Lemma~\ref{lem:colour-words}, and obtain a formula of the same shape
as~\eqref{eq:grouped} \cite[Corollary~4.2]{PdAW}, with an additional parameter attached to the
crossings between different components. In comparison, Theorems~\ref{thm:state-sum}
and~\ref{thm:realization} show that state sums of this type govern all diagonal enhancements of
charge-conserving Yang--Baxter operators, with independent HOMFLY--PT specializations on the
different colours and arbitrary couplings $\omega_{st}$, and that each of them occurs. On the
other hand, the specializations $P_L(A,z)$ that occur satisfy $A\in\pm Q^{\Z}$ for
$z=Q-Q^{-1}$, whereas Markov traces on Hecke algebras~\cite{Jon} and on Yokonuma--Hecke algebras
produce arbitrary values of $A$.

\subsection{Further remarks}
\label{subsec:further}
\emph{The case $n=0$.} When $n=0$ the connected invariant vanishes identically, as does the
unreduced invariant attached to $\mathfrak{gl}_{1|1}$, whose reduced version is the
Alexander--Conway polynomial \cite{KS,QS}. Aizawa et al.\ obtained the Alexander polynomial
from $R_{2.2}$, for which $n=0$, by means of a state model \cite[Section~4.2]{AHKO}. Reduced
versions of the invariants of Theorem~\ref{thm:state-sum}, obtained by cutting open one strand,
are a natural next step.

\emph{Non-diagonal and non-generic enhancements.} A charge-conserving operator can admit
enhancements whose $\mu$ is not diagonal. For instance, if $R=\lambda\,\id$, then
$(R,\mu,\lambda,\beta)$ is an enhanced Yang--Baxter operator for every $\mu\in\GL(V)$ with
$\tr(\mu)=\beta$. By the discussion leading to~\eqref{eq:cc-intro}, operators that are not
charge-conserving with respect to the eigenlines of $\mu$ can only appear in enhancements whose
$\mu$ is not diagonalizable or has non-generic eigenvalues. This happens for the
two-dimensional solution $R_{1.2}$ of \cite[Section~3]{AHKO}, which depends on a parameter
$q\in\C^\times$ and satisfies $R_{1.2}(e_+\otimes e_+)=e_+\otimes e_++e_-\otimes e_-$. For every
$\beta\in\C^\times$ it admits the enhancement $(R_{1.2},\mu,\sqrt q,\beta)$ with
$\mu=\operatorname{diag}(\beta/\sqrt q,-\beta/\sqrt q)$ \cite[Table~1]{AHKO}, whose eigenvalues
satisfy $\mu_+^2=\mu_-^2$, although $R_{1.2}$ is not charge-conserving with respect to the
eigenlines $\C e_+$ and $\C e_-$ of $\mu$. The classification of these enhancements remains
open.

\section*{Acknowledgements}
We used Claude Opus 5.5 during the development and drafting of this paper, to implement some
proof strategies and to explore some computational examples. All statements and proofs were
checked by the authors, who are responsible for them. C.G.\ was partially supported by Grant
INV-2025-213-3452 from the School of Science of Universidad de los Andes.

\bibliographystyle{amsplain}
\bibliography{references}

\end{document}